\documentclass{svjour3}

\usepackage{dpr_fec}
\usepackage{longtable}
\usepackage[table]{xcolor}
\usepackage{caption}
\usepackage{subcaption}
\usepackage{wrapfig}
\hypersetup{colorlinks}
\smartqed
\allowdisplaybreaks

\begin{document}

%*******
% Title
%*******
\title{Exploiting Negative Curvature in Deterministic and Stochastic Optimization}
\author{Frank E.~Curtis
        \and
        Daniel P.~Robinson}
\institute{F.~E.~Curtis \at 
           Department of Industrial and Systems Engineering \\
           Lehigh University \\
           \email{frank.e.curtis@gmail.com} \\
           This author was supported in part by the U.S.~Department of Energy under Grant No.~DE-SC0010615 and by the U.S.~National Science Foundation under Grant No.~CCF-1618717.
           \and
           D.~P.~Robinson \at
           Department of Applied Mathematics and Statistics \\
           Johns Hopkins University \\
           \email{daniel.p.robinson@gmail.com}}
\date{\today}
\maketitle

%**********
% Abstract
%**********
\begin{abstract}
  This paper addresses the question of whether it can be beneficial for an optimization algorithm to follow directions of negative curvature.  Although  prior work has established convergence results for algorithms that integrate both descent and negative curvature steps, there has not yet been extensive numerical evidence showing that such methods offer consistent performance improvements.  In this paper, we present new frameworks for combining descent and negative curvature directions: alternating two-step approaches and dynamic step approaches.  The aspect that distinguishes our approaches from ones previously proposed is that they make algorithmic decisions based on (estimated) upper-bounding models of the objective function.  A consequence of this aspect is that our frameworks can, in theory, employ fixed stepsizes, which makes the methods readily translatable from deterministic to stochastic settings.  For deterministic problems, we show that instances of our dynamic framework yield gains in performance compared to related methods that only follow descent steps.  We also show that gains can be made in a stochastic setting in cases when a standard stochastic-gradient-type method might make slow progress.
  \keywords{nonconvex optimization \and
            second-order methods \and
            modified Newton methods \and
            negative curvature \and
            stochastic optimization \and            
            machine learning}
  \subclass{49M05 \and
            49M15 \and
            49M37 \and
            65K05 \and
            90C15 \and
            90C26 \and
            90C30 \and
            90C53}
\end{abstract}

%*********
% Section
%*********
\section{Introduction}\label{sec.introduction}

There has been a recent surge of interest in solving nonconvex optimization problems.  A prime example is the dramatic increase in interest in the training of deep neural networks.
%, a subject that admits very challenging nonconvex optimization problems.  
Another example is the task of clustering data that arise from the union of low dimensional subspaces.  In this setting, the nonconvexity typically results from sophisticated modeling approaches that attempt to accurately capture corruptions in the data~\cite{elhamifar2009sparse,JiaRV17}.  It is now widely accepted that the design of new methods for solving nonconvex problems (at least locally) is sorely needed.

First consider deterministic optimization problems.  For solving such problems, most algorithms designed for minimizing smooth objective functions only ensure convergence to first-order stationarity, i.e., that the gradient of the objective asymptotically vanishes.  This characterization is certainly accurate for line search methods, which seek to reduce the objective function by searching along descent directions.  Relatively few researchers have designed line search (or other, such as trust region or regularization) algorithms that generate iterates that provably converge to second-order stationarity.  The reason for this is three-fold: $(i)$~such methods are more complicated and expensive, necessarily involving the computation of directions of negative curvature when they exist; $(ii)$ methods designed only to achieve first-order stationarity rarely get stuck at saddle points that are first-order, but not second-order stationary \cite{pmlr-v49-lee16}; and $(iii)$ there has not been sufficient evidence showing benefits of integrating directions of negative curvature.

For solving stochastic optimization problems, the methods most commonly invoked are variants of the stochastic gradient (SG) method. During each iteration of SG, a stochastic gradient is computed and a step opposite that direction is taken to obtain the next iterate.  Even for nonconvex problems, convergence guarantees (e.g., in expectation or almost surely) for SG methods to first-order stationarity have been established under reasonable assumptions; e.g., see~\cite{BottCurtNoce16}.  In fact, SG and its variants represent the current state-of-the-art for training deep neural networks.  As for methods that compute and follow negative curvature directions, it is no surprise that such methods (e.g.,~\cite{liu2017noisy}) have not been used or studied extensively since practical benefits have not even been shown in deterministic optimization.

The main purpose of this paper is to revisit and provide new perspectives on the use of negative curvature directions in deterministic and stochastic optimization.  Whereas previous work in deterministic settings has focused on line search and other %et~al.~
methods, we focus on methods that attempt to construct \emph{upper-bounding models} of the objective function.  This allows them, at least in theory, to employ \emph{fixed stepsizes} while ensuring convergence guarantees.  For a few instances of such methods of interest, we provide theoretical convergence guarantees and empirical evidence showing that an optimization process can benefit by following negative curvature directions.  The fact that our methods might employ fixed stepsizes is also important as it allows us to offer new strategies for stochastic optimization where, e.g., line search strategies are not often viable.  
%Overall, the theme of this paper can be summarized as the following: \textit{In contrast to methods designed to avoid ill-effects associated with the presence of negative curvature, our methods {\rm exploit} directions of negative curvature to achieve gains in performance.}

%************
% Subsection
%************
\subsection{Contributions}\label{sec.contributions}

The contributions of our work are the following.

\bitemize
  \item For deterministic optimization, we first provide conditions on descent and negative curvature directions that allow us to guarantee convergence to second-order stationarity with a simple two-step iteration with fixed stepsizes; see~\S\ref{sec.deterministic.two-step}.  
  Using 
  %\item Our second method for deterministic optimization uses the
  the two-step method as motivation, we then propose a dynamic choice for the direction and stepsize; see~\S\ref{sec.deterministic.dynamic}.  In particular, the dynamic algorithm makes decisions based on which available step appears to offer a more significant objective reduction. The details of our dynamic strategy represent the main novelty of this paper with respect to deterministic optimization.
  \item We prove convergence rate guarantees for our deterministic optimization methods that provide upper bounds for the numbers of iterations required to achieve (approximate) first- and second-order stationarity; see~\S\ref{sec.deterministic.complexity}.  We follow this with a discussion, in \S\ref{sec.deterministic.discussion}, on the issue of how methods might behave in the neighborhood of so-called strict saddle points, a topic of much interest in the literature.
  \item In \S\ref{sec.deterministic.step_computation} and \S\ref{sec.deterministic.numerical}, we discuss different techniques for computing the search directions in our deterministic algorithms and provide the results of numerical experiments showing the benefits of following negative curvature directions, as opposed to only descent directions.
  \item For solving stochastic optimization problems, we propose two methods.  Our first method shows that one can maintain the convergence guarantees of a stochastic gradient method by adding an appropriately scaled negative curvature direction for a stochastic Hessian estimate; see~\S\ref{sec.stochastic.two-step}.  This approach can be seen as a refinement of that in~\cite{NeelVilnLeSutsKaisKuraMart15}, which adds noise to each SG step.
  \item Our second method for stochastic optimization is an adaptation of our dynamic (deterministic) method when stochastic gradient and Hessian estimates are involved; see~\S\ref{sec.stochastic.dynamic}.  Although we are unable to establish a convergence theory for this approach, we do illustrate some gain in performance in neural network training; see~\S\ref{sec.stochastic.numerical}.  We view this as a first step in the design of a practical algorithm for stochastic optimization that efficiently exploits negative curvature.
\eitemize

Computing directions of negative curvature carries an added cost that should be taken into consideration when comparing algorithm performance.  We remark along with the results of our numerical experiments how these added costs can be worthwhile and/or marginal relative to the other per-iteration costs.% in an algorithm.  It is also worthwhile to mention that some researchers argue that negative curvature in, say, the training of deep networks should not be explored.  They promote this view since some nonlinearity is introduced by the choice of activation function, and is not a true feature of the underlying problem of interest.  We disagree.  After all, regardless of its origin, nonlinearity is a part of the optimization problem to be solved, and negative curvature directions can offer the possibility of reaching in productive directions that might not otherwise be explored.~\footnote{dpr:this last sentence is a bit off. Also, in general, should we bring this up?  Can we be blamed again for trying to shoot the complan down before it is made.}

%************
% Subsection
%************
\subsection{Prior Related Work}\label{sec.prior}

For solving deterministic optimization problems, relatively little research has been directed towards the use of negative curvature directions.  Perhaps the first exception is the work in~\cite{more1979use} in which convergence to second-order stationary points is proved using a curvilinear search formed by descent and negative curvature directions.  In a similar vein, the work in~\cite{forsgren1995computing} offers similar convergence properties based on using a curvilinear search, although the primary focus in that work is describing how a partial Cholesky factorization of the Hessian matrix could be used to compute descent and negative curvature directions.  More recently, a linear combination of descent and negative curvature directions was used in an optimization framework to establish convergence to second-order stationary solutions under loose assumptions~\cite{GilKR16,GilKR16c}, and a strategy for how to use descent and negative curvature directions was combined with a backtracking linesearch to provide worst case iteration bounds in~\cite{royer2017complexity}. 
For further instances of work employing negative curvature directions, see \cite{Birgin2001,goldfarb1980curvilinear,doi:10.1080/10556780008805794,doi:10.1137/15M1022100}.  Importantly, none of the papers above (or any others to the best of our knowledge) have established consistent gains in computational performance as a result of using negative curvature directions.

Another recent trend in the design of deterministic methods for nonconvex optimization is to focus on the ability of an algorithm to escape regions around a saddle point, i.e., a first-order stationary point that is not a minimizer nor a maximizer.  A prime example of this trend is the work in \cite{PateMokhRibe17} in which a standard type of regularized Newton method is considered.  (For a more general presentation of regularized Newton methods of which the method in \cite{PateMokhRibe17} is a special case, see, e.g.,  \cite{NoceWrig06}.)  While the authors do show a probabilistic convergence result for their method to (approximate) second-order stationarity, their main emphasis is on the number of iterations required to escape neighborhoods of saddle points.  This result, similar to those in some (but not all) recent papers discussing the behavior of descent methods when solving nonconvex problems (e.g., see \cite{jin2017escape,pmlr-v49-lee16}), requires that all saddle points are \emph{non-degenerate} in the sense that the Hessian of the objective at any such point does not have any zero eigenvalues.  In particular, in~\cite{jin2017escape} they show that a perturbed form of gradient descent converges to a second-order stationary point in a number iterations that depends poly-logarithmically on the dimension of the problem.  Our convergence theory does not require such a non-degeneracy assumption; see~\S\ref{sec.deterministic.discussion} for additional discussion comparing our results to others. (See also \cite{LeePanPilSimJorRec17}, where it is shown, without a non-degeneracy assumption and for almost all starting points, that certain first-order methods do not converge to saddle points at which the Hessian of the objective has a negative eigenvalue.  That said, this work does not prove bounds on the number of iterations such a method might spend near saddle points.)

For solving stochastic optimization problems, there has been little work that focuses explicitly on the use of directions of negative curvature.  For two examples, see \cite{NIPS2015_5870,liu2017noisy}.  Meanwhile, it has recently become clear that the nonconvex optimization problems used to train deep neural networks have a rich and interesting landscape~\cite{dauphin2014identifying,keskar2016large}.  Instead of using negative curvature to their advantage, modern methods ignore it, introduce random perturbations (e.g., see \cite{GeHuanJinYuan15}), or employ regularized/modified Newton methods that attempt to avoid its potential ill-effects (e.g., \cite{dauphin2014identifying,martens2010deep}).  Although such approaches are reasonable, one might intuitively expect better performance if directions of negative curvature are exploited instead of avoided. In this critical respect, the methods that we propose are different from these prior algorithms, with the exception of~\cite{liu2017noisy}.

%*********
% Section
%*********
\section{Deterministic Optimization}\label{sec.deterministic}

Consider the unconstrained optimization problem
\bequation\label{prob.nlp}
  \min_{x\in\R{n}}\ f(x),
\eequation
where the objective function $f : \R{n} \to \R{}$ is twice continuously differentiable and bounded below by a scalar $f_{\inf} \in \R{}$.  We define the gradient function $g := \nabla f$ and Hessian function $H := \nabla^2 f$.  We assume that both of these functions are Lipschitz continuous on the path defined by the iterates computed in an algorithm, the gradient function $g$ with Lipschitz constant $L \in (0,\infty)$ and the Hessian function $H$ with Lipschitz constant $\sigma \in (0,\infty)$.  Given an invertible matrix $M \in \R{n\times n}$, its Euclidean norm condition number is written as $\kappa(M) = \|M\|_2\|M^{-1}\|_2$.  Given a scalar $\lambda \in \R{}$, we define $(\lambda)_- := \min\{0,\lambda\}$.

In the remainder of this section, we present a two-step method that is guaranteed to converge toward second-order stationarity, as well as a dynamic approach that chooses between two types of steps at each point based on lower bounds on objective function decrease.  Both algorithms are presented in a generic manner that offers flexibility in the ways the steps are computed.

%************
% Subsection
%************
\subsection{Two-Step Method}\label{sec.deterministic.two-step}

Our first method alternates between negative curvature and descent steps using fixed stepsizes for each. (Either can be taken first; arbitrarily, we state our algorithm and analysis assuming that one starts with a negative curvature step.)  At a given iterate $x_k \in \R{n}$, let~$\lambda_k$ denote the left-most eigenvalue of $H(x_k)$.  If $\lambda_k \geq 0$ (i.e., $H(x_k) \succeq 0$), the algorithm sets $d_k \gets 0$; otherwise, $d_k$ is computed so that
\bsubequations\label{def.d}
  \begin{align}
    d_k^TH(x_k)d_k &\leq \gamma\lambda_k\|d_k\|_2^2 < 0, \label{def.d.curv} \\
    g(x_k)^Td_k &\leq 0, \label{def.d.descent} \\
    \text{and}\ \ \|d_k\|_2 &\leq \theta |\lambda_k|,\label{def.d.norm}
  \end{align}
\esubequations
for some $\gamma \in (0,1]$ and $\theta \in (0,\infty)$ that are independent of $k$.  A step in this direction is then taken to obtain $\xhat_k \gets x_k + \beta d_k$ for some $\beta \in(0,(3\gamma)/(\sigma\theta))$.  At this point, if $g(\xhat_k) = 0$, then the algorithm sets $\shat_k \gets 0$; otherwise, $\shat_k$ is computed satisfying
\bequation\label{def.s}
  \frac{-g(\xhat_k)^T\shat_k}{\|\shat_k\|_2\|g(\xhat_k)\|_2} \geq \delta\ \ \text{and}\ \ \zeta \leq \frac{\|\shat_k\|_2}{\|g(\xhat_k)\|_2} \leq \eta
\eequation
for some $(\delta,\zeta) \in (0,1] \times (0,1]$ and $\eta \in [1,\infty)$ that are independent of $k$.  The iteration ends by taking a step along this direction to obtain $x_{k+1} \gets \xhat_k + \alpha \shat_k \equiv x_k + \alpha\shat_k + \beta d_k$ for some $\alpha \in(0,(2\delta\zeta)/(L\eta^2))$.  We remark that \eqref{def.d} and~\eqref{def.s} are satisfiable; e.g., if $\shat_k = -g(\xhat_k)$ and $d_k$ is an eigenvector corresponding to the left-most eigenvalue $\lambda_k$ scaled so that $g(x_k)^Td_k \leq 0$ and $\|d_k\|_2 = |\lambda_k|$, then \eqref{def.d} and~\eqref{def.s} are satisfied with $\gamma = \theta = \delta = \zeta = \eta = 1$.  For further remarks on step computation techniques, see \S\ref{sec.deterministic.step_computation}.

\balgorithm[ht]
  \caption{Two-Step Method}
  \label{alg.2step}
  \balgorithmic[1]
    \Require $x_1 \in \R{n}$, $\alpha\in(0, (2\delta\zeta)/(L\eta^2))$, and $\beta\in(0,(3\gamma)/(\sigma\theta))$
    \For{\textbf{all} $k \in \{1,2,\dots\} =: \N{}_+$}
        \State \textbf{if} $\lambda_k \geq 0$ \textbf{then} set $d_k \gets 0$ \textbf{else} set $d_k$ satisfying~\eqref{def.d}
        \State set $\xhat_k \gets x_k + \beta d_k$
        \State \textbf{if} $g(\xhat_k) = 0$ \textbf{then} set $\shat_k \gets 0$ \textbf{else} set $\shat_k$ satisfying~\eqref{def.s}
        \State \textbf{if} $d_k = \shat_k = 0$ \textbf{then return} $x_k$
        \State set $x_{k+1} \gets \xhat_k + \alpha_k \shat_k \equiv x_k + \alpha \shat_k + \beta d_k$ \label{step.update}
     \EndFor  
  \ealgorithmic
\ealgorithm

We now show that Algorithm~\ref{alg.2step} converges toward second-order stationarity.  Critical for this analysis are bounds on $(\alpha,\beta)$ that we have stated (and are also stated in Algorithm~\ref{alg.2step}).  Also, for the analysis, it is convenient to define
\bequationNN
  \Dcal := \{k \in \N{}_+ : d_k \neq 0\} \equiv \{k \in \N{}_+ : \lambda_k < 0\}
\eequationNN
along with the indicator $\Ical_\Dcal(k)$, which evaluates as 1 if $k \in \Dcal$ and as 0 otherwise.

\btheorem\label{th.deterministic}
  If Algorithm~\ref{alg.2step} terminates finitely in iteration $k \in \N{}_+$, then $g(x_k) = 0$ and $\lambda_k \geq 0$, i.e., $x_k$ is second-order stationary.  Otherwise, the computed iterates satisfy
  \bequation\label{eq.limits}
    \lim_{k\to\infty} \|g(x_k)\|_2 = 0\ \ \text{and}\ \ \liminf_{k\to\infty} \lambda_k \geq 0.
  \eequation
\etheorem
\bproof
  Algorithm~\ref{alg.2step} terminates finitely if and only if, for some $k \in \N{}_+$, $d_k = \shat_k = 0$.   We can then observe that $d_k = 0$ if and only if $\lambda_k \geq 0$, while $\shat_k = 0$ if and only if $0 = g(\xhat_k) = g(x_k)$, where the last equality holds because $\xhat_k = x_k$ whenever $d_k = 0$.
%  This can only occur if $\lambda_k \geq 0$ and $g(x_k) = 0$, where we have used the fact that when $d_k = 0$ it follows that $\xhat_k = x_k$.  
These represent the desired conclusions for this case.  Otherwise, if Algorithm~\ref{alg.2step} does not terminate finitely, consider arbitrary $k \in \N{}_+$.  If $k\notin\Dcal$, then $d_k = 0$ and $\xhat_k = x_k$, meaning that $f(\xhat_k) = f(x_k)$.  Otherwise, if $k\in\Dcal$, then $d_k \neq 0$ and $\lambda_k < 0$, and, by the step computation conditions in \eqref{def.d}, it follows that
  \bequationNN
    \baligned
    f(\xhat_k)
      &\leq f(x_k) + g(x_k)^T(\beta d_k) + \thalf (\beta d_k)^TH(x_k)(\beta d_k) + \tfrac16 \sigma \|\beta d_k\|_2^3 \\
      &\leq f(x_k) + \thalf \beta^2 \gamma \lambda_k \|d_k\|_2^2 + \tfrac16 \sigma \beta^3 \theta^3 |\lambda_k|^3 \\
      &\leq f(x_k) - \thalf\beta^2\theta^2 \big(\gamma - \tfrac13 \sigma\beta \theta \big)  |\lambda_k|^3 \\
      &= f(x_k) - c_1(\beta) |\lambda_k|^3,
    \ealigned
  \eequationNN
  where $c_1(\beta) : = \thalf\beta^2\theta^2 \big( \gamma - \tfrac13 \sigma\beta \theta \big) \in (0,\infty)$.  Similarly, it follows from~\eqref{def.s} that
  \bequationNN
    \baligned
    f(x_{k+1})
      &\leq f(\xhat_k) + g(\xhat_k)^T(\alpha\shat_k) + \thalf L \|\alpha\shat_k\|_2^2 \\
      &\leq f(\xhat_k) - \alpha \delta \zeta \|g(\xhat_k)\|_2^2  + \thalf L \alpha^2 \eta^2 \|g(\xhat_k)\|_2^2 \\
      &= f(\xhat_k) - \alpha (\delta\zeta - \thalf L \alpha \eta^2) \|g(\xhat_k)\|_2^2 \\
      &= f(\xhat_k) - c_2(\alpha)  \|g(\xhat_k)\|_2^2,
    \ealigned
  \eequationNN
  where $c_2(\alpha) := \alpha (\delta\zeta - \thalf L \alpha \eta^2) \in (0,\infty)$.   Overall,
  \bequation\label{eq.twostep_reduction}
    f(x_{k+1}) \leq f(x_k) - \Ical_\Dcal(k) c_1(\beta) |\lambda_k|^3 - c_2(\alpha) \|g(\xhat_k)\|_2^2.
  \eequation
  Now observe that, for any $\ell \in \N{}_+$, it follows that
  \bequationNN
    \baligned
      f(x_1) - f(x_{\ell+1}) &= \sum_{k=1}^\ell (f(x_k) - f(x_{k+1})) \\
      &\geq \sum_{k=1}^\ell \Big( \Ical_\Dcal(k) c_1(\beta) |\lambda_k|^3 + c_2(\alpha)\|g(\xhat_k)\|_2^2 \Big).
    \ealigned
  \eequationNN
  This inequality and $f$ being bounded below show that
  \bsubequations
    \begin{align}
      \sum_{k=1, k\in\Dcal}^\infty |\lambda_k|^3 \equiv \sum_{k=1}^\infty |(\lambda_k)_-|^3 &< \infty \label{eq.sums.lambda} \\ 
      \text{and} \ \ \sum_{k=1}^\infty \|g(\xhat_k)\|_2^2 &< \infty. \label{eq.sums.g}
    \end{align}
  \esubequations
  The latter bound yields
  \bequation\label{eq.grad-xhat-0}
    \lim_{k\to\infty} \|g(\xhat_k)\|_2 = 0,
  \eequation
  while the former bound and~\eqref{def.d.norm} yield
  \bequationNN
    \sum_{k=1}^\infty \|\xhat_k - x_k\|_2^3 = \beta^3 \sum_{k=1}^\infty\|d_k\|_2^3 
    \leq \beta^3\theta^3 \sum_{k=1}^\infty |(\lambda_k)_-|^3 < \infty,
  \eequationNN
  from which it follows that
  \bequation\label{eq.xhat-to-x}
    \lim_{k\to\infty} \|\xhat_k - x_k\|_2 = 0.
  \eequation
  It follows from Lipschitz continuity of $g$ along with \eqref{eq.grad-xhat-0} and \eqref{eq.xhat-to-x} that
  \bequationNN
    \baligned
    0 
      &\leq \limsup_{k\to\infty} \|g(x_k)\|_2 \\
      &=    \limsup_{k\to\infty} \|g(x_k) - g(\xhat_k) + g(\xhat_k)\|_2 \\
      &\leq \limsup_{k\to\infty} \|g(x_k) - g(\xhat_k)\|_2 + \limsup_{k\to\infty} \|g(\xhat_k)\|_2 \\
      &\leq L \limsup_{k\to\infty} \|x_k - \xhat_k\|_2 + \limsup_{k\to\infty} \|g(\xhat_k)\|_2 = 0,
    \ealigned
  \eequationNN
  which implies the first limit in \eqref{eq.limits}.  Finally, in order to derive a contradiction, suppose that $\liminf_{k\to\infty} \lambda_k < 0$, meaning that there exists some $\epsilon > 0$ and infinite index set $\Kcal \subseteq \Dcal$ such that $\lambda_k \leq -\epsilon$ for all $k \in \Kcal$.  This implies that
  \bequationNN
    \sum_{k=1}^\infty |(\lambda_k)_-|^3 
    \geq \sum_{k\in\Kcal} |(\lambda_k)_-|^3 
    \geq \sum_{k\in\Kcal} \epsilon^3 = \infty,
  \eequationNN
  contradicting \eqref{eq.sums.lambda}.  This yields the second limit in \eqref{eq.limits}. \qed
\eproof

There are two potential weaknesses of this two-step approach.  First, it simply alternates back-and-forth between descent and negative curvature directions, which might not always lead to the most productive step from each point.  Second, even though our analysis holds for all stepsizes $\alpha$ and $\beta$ in the intervals provided in Algorithm~\ref{alg.2step}, the algorithm might suffer from poor performance if these values are chosen poorly.  We next present a method that addresses these weaknesses.

%************
% Subsection
%************
\subsection{Dynamic Method}\label{sec.deterministic.dynamic}

Suppose that, in any iteration $k \in \N{}_+$ when $\lambda_k < 0$, one computes a nonzero direction of negative curvature satisfying \eqref{def.d.curv}--\eqref{def.d.descent} for some $\gamma \in (0,1]$.  Suppose also that, if $g(x_k) \neq 0$, then one computes a nonzero direction~$s_k$ satisfying the equivalent of the first condition in \eqref{def.s}, namely, for some $\delta \in (0,1]$,
\bequation\label{def.s.adapt}
  -g(x_k)^Ts_k \geq \delta \|s_k\|_2\|g(x_k)\|_2.
\eequation
Defining, for $(L_k,\sigma_k) \in (0,\infty) \times (0,\infty)$, the model reductions
\bequationNN
  \baligned
  m_{s,k}(\alpha) &:= -\alpha g(x_k)^T s_k - \thalf L_k \alpha^2 \|s_k\|_2^2 \\ \text{and}\ \  
  m_{d,k}(\beta)  &:= -\beta  g(x_k)^T d_k - \thalf \beta^2 d_k^T H(x_k) d_k - \tfrac{\sigma_k}{6} \beta^3 \|d_k\|_2^3,
  \ealigned
\eequationNN
ones finds that, if $L_k \geq L$ and $\sigma_k \geq \sigma$, then
\bsubequations\label{mod.dec.alpha_beta}
  \begin{align}
    f(x_k + \alpha s_k) &\leq f(x_k) - m_{s,k}(\alpha) \label{mod.dec.alpha} \\ \text{and}\ \ 
    f(x_k + \beta d_k)  &\leq f(x_k) - m_{d,k}(\beta). \label{mod.dec.beta}
  \end{align}
\esubequations
These two inequalities suggest that, during iteration $k \in \N{}_+$, one could choose which of the two steps ($s_k$ or $d_k$) to take based on which model reduction predicts the larger decrease in the objective.  One can verify that the reductions $m_{s,k}$ and $m_{d,k}$ are maximized (over the positive real numbers) by
\bequation\label{eq.alpha.beta}
  \alpha_k := \frac{- g(x_k)^T s_k}{L_k\|s_k\|_2^2} \ \ \text{and}\ \ \beta_k := \frac{\left(-c_k + \sqrt{c_k^2 - 2\sigma_k\|d_k\|_2^3 g(x_k)^T d_k}\right)}{\sigma_k\|d_k\|_2^3},
\eequation
where $c_k := d_k^T H(x_k)d_k$ is a measure of $H(x_k)$-curvature for $d_k$.

Algorithm~\ref{alg.dynamic}, stated below, follows this dynamic strategy of choosing between~$s_k$ and $d_k$ for all $k \in \N{}_+$.  It also involves dynamic updates for Lipschitz constant estimates, represented in iteration $k \in \N{}_+$ by $L_k$ and $\sigma_k$.  In this deterministic setting, a step is only taken if it yields an objective function decrease.  Otherwise, a null step is effectively taken and a Lipschitz constant estimate is increased.

\balgorithm[ht]
  \caption{Dynamic Method}
  \label{alg.dynamic}
  \balgorithmic[1]
    \Require $x_1 \in \R{n}$ and $(\rho,L_1,\sigma_1) \in (1,\infty) \times (0,\infty) \times (0,\infty)$
    \For{\textbf{all} $k \in \N{}_+$}
        \State \textbf{if} $\lambda_k \geq 0$ \textbf{then} set $d_k \gets 0$ \textbf{else} set $d_k$ satisfying~\eqref{def.d.curv}--\eqref{def.d.descent}
        \State \textbf{if} $g(x_k) = 0$ \textbf{then} set $s_k \gets 0$ \textbf{else} set $s_k$ satisfying~\eqref{def.s.adapt}
        \State \textbf{if} $d_k = s_k = 0$ \textbf{then return} $x_k$\label{step.return}
        \Loop\label{step.loop}
        \State compute $\alpha_k > 0$ and $\beta_k > 0$ from~\eqref{eq.alpha.beta}\label{step.alphabeta}
        \If{$m_{s,k}(\alpha_{k}) \geq m_{d,k}(\beta_{k})$}
           \If{\eqref{mod.dec.alpha} holds with $\alpha = \alpha_k$}
              \State set $x_{k+1} \gets x_k + \alpha_k s_k$
              \State \textbf{exit loop}
           \Else
              \State set $L_k \gets \rho L_k$ \label{step.L-increase}
           \EndIf
        \Else
            \If{\eqref{mod.dec.beta} holds with $\beta = \beta_k$}
              \State set $x_{k+1} \gets x_k + \beta_k d_k$
              \State \textbf{exit loop}
           \Else
              \State set $\sigma_k \gets \rho \sigma_k$ \label{step.sigma-increase}
           \EndIf
        \EndIf
        \EndLoop
        \State set $(L_{k+1},\sigma_{k+1}) \in (0,L_k] \times (0,\sigma_k]$ \label{step.set_L_sigma}
     \EndFor
  \ealgorithmic
\ealgorithm

In the next two results, we establish that Algorithm~\ref{alg.dynamic} is well-defined and that it has convergence guarantees on par with Algorithm~\ref{alg.2step} (recall Theorem~\ref{th.deterministic}).

\blemma
  Algorithm~\ref{alg.dynamic} is well defined in the sense that it either terminates finitely or generates infinitely many iterates.  In addition, at the end of each iteration $k \in \N{}_+$,
  \bequation\label{eq.L_sigma_upper}
    L_k \leq L_{\max} := \max\{L_1,\rho L\}\ \ \text{and}\ \ \sigma_k \leq \sigma_{\max} := \max\{\sigma_1,\rho \sigma\}.
  \eequation
\elemma
\bproof
We begin by showing that the \textbf{loop} in Step~\ref{step.loop} terminates finitely anytime it is entered. For a proof by contradiction, if
% If it does not, then it enters the \textbf{loop}.  If
that loop were never to terminate, then the updates to $L_k$ and/or~$\sigma_k$ would cause at least one of them to become arbitrarily large.  Since~\eqref{mod.dec.alpha} holds whenever $L_k \geq L$ and~\eqref{mod.dec.beta} holds whenever $\sigma_k \geq \sigma$, it follows that the loop would eventually terminate, thus reaching a contradiction.  Therefore, we have proved that the \textbf{loop} in Step~\ref{step.loop} terminates finitely anytime it is entered, and moreover that~\eqref{eq.L_sigma_upper} holds. This completes the proof once we observe that during iteration $k \in \N{}_+$, Algorithm~\ref{alg.dynamic} might finitely terminate in Step~\ref{step.return}.  
 \qed
\eproof

\btheorem\label{th.deterministic.adaptive}
  If Algorithm~\ref{alg.dynamic} terminates finitely in iteration $k \in \N{}_+$, then $g(x_k) = 0$ and $\lambda_k \geq 0$, i.e., $x_k$ is second-order stationary.  Otherwise, the computed iterates satisfy
  \bequation\label{eq.limits.dynamic}
    \lim_{k\to\infty} \|g(x_k)\|_2 = 0\ \ \text{and}\ \ \liminf_{k\to\infty} \lambda_k \geq 0.
  \eequation
\etheorem
\bproof
  Algorithm~\ref{alg.dynamic} terminates finitely only if, for some $k \in \N{}_+$, $d_k = s_k = 0$.  This can only occur if $\lambda_k \geq 0$ and $g(x_k) = 0$, which are the desired conclusions.  Otherwise, Algorithm~\ref{alg.dynamic} does not terminate finitely and during each iteration at least one of $s_k$ and $d_k$ is nonzero.  We use this fact in the following arguments.
 
  Consider arbitrary $k \in \N{}_+$.  If $s_k \neq 0$, then the definition of~$\alpha_k$ in \eqref{eq.alpha.beta} and the step computation condition on $s_k$ in~\eqref{def.s.adapt} ensure that
  \bequation\label{bd.dec.alpha}
    m_{s,k}(\alpha_k) = \frac{1}{2L_k} \left(\frac{g(x_k)^T s_k}{\|s_k\|_2}\right)^2 \geq \frac{\delta^2}{2L_k} \|g(x_k)\|_2^2.
  \eequation
  Similarly, if $d_k \neq 0$, then by the fact that $\beta_k$ maximizes $m_{d,k}(\beta)$ over $\beta > 0$, \eqref{def.d.curv}--\eqref{def.d.descent}, and defining $\hat\beta_k := -2d_k^T H(x_k)d_k/(\sigma_k\|d_k\|_2^3) > 0$, one finds
  \begin{align}
    m_{d,k}(\beta_{k})
      &\geq m_{d,k}(\hat\beta_k) \nonumber \\
      &\geq -\thalf\hat\beta_k^2 d_k^T H(x_k) d_k - \tfrac16 \sigma_k \hat\beta_k^3 \|d_k\|_2^3 \nonumber \\
      &=    -\frac{2(d_k^T H(x_k) d_k)^3}{3\sigma_k^2\|d_k\|_2^6}
       \geq  \frac{2\gamma^3}{3\sigma_k^2} |\lambda_k|^3
       = \frac{2\gamma^3}{3\sigma_k^2} |(\lambda_k)_-|^3. \label{bd.dec.beta}
  \end{align}
%  One can extend this inequality to all cases (i.e., not only when $d_k \neq 0$) using the fact that $()_-$ whenver $d_k = 0$, via
%  \bequation\label{bd.dec.beta}
%    m_{d,k}(\beta_{k}) \geq \frac{2\gamma^3}{3\sigma_k^2} |(\lambda_k)_-|^3\ \ \text{for all}\ \ k \in \N{}.
 % \eequation
  Overall, for all $k \in \N{}_+$, since at least one of $s_k$ and $d_k$ is nonzero, $\|g_k\|_2 = 0$ if and only if $\|s_k\|_2 = 0$, and $|(\lambda_k)_-| = 0$ if and only if $\|d_k\|_2 = 0$, it follows that
  \bequation\label{dec.dynamic}
      f(x_k) - f(x_{k+1}) \geq \max\left\{\frac{\delta^2}{2L_k} \|g(x_k)\|_2^2, \frac{2\gamma^3}{3\sigma_k^2} |(\lambda_k)_-|^3 \right\}.
  \eequation
  Indeed, to show that~\eqref{dec.dynamic} holds, let us consider two cases.  First, suppose that the update $x_{k+1} \gets x_k + \alpha_k s_k$ is completed, meaning that~\eqref{mod.dec.alpha} holds with $\alpha = \alpha_k$ and $m_{s,k}(\alpha_k) \geq m_{d,k}(\beta_k)$.  Combining these facts with~\eqref{bd.dec.alpha} and~\eqref{bd.dec.beta} establishes~\eqref{dec.dynamic} in this case.  Second, suppose that $x_{k+1} \gets x_k + \beta_k d_k$ is completed, meaning that~\eqref{mod.dec.beta} holds with $\beta = \beta_k$ and $m_{s,k}(\alpha_k) < m_{d,k}(\beta_k)$. Combining these facts with~\eqref{bd.dec.alpha} and~\eqref{bd.dec.beta} establishes~\eqref{dec.dynamic} in this case. Thus, \eqref{dec.dynamic} holds for all $k \in \N{}_+$.

  It now follows from~\eqref{dec.dynamic}, the bounds in \eqref{eq.L_sigma_upper}, and a proof similar to that used in Theorem~\ref{th.deterministic} (in particular, to establish~\eqref{eq.sums.lambda} and~\eqref{eq.sums.g}) that
  \bequationNN
    \sum_{k=1}^\infty \|g(x_k)\|_2^2 < \infty\ \ \text{and} \ \ \sum_{k=1}^\infty |(\lambda_k)_-| < \infty. 
  \eequationNN
  One may now establish the desired results in~\eqref{eq.limits.dynamic} using the same arguments as used in the proof of Theorem~\ref{th.deterministic}. \qed
\eproof

Let us add a few remarks about Algorithm~\ref{alg.dynamic}.
\bitemize
  \item Our convergence theory allows $L_{k+1} \gets L_k$ and $\sigma_{k+1} \gets \sigma_k$ in Step~\ref{step.set_L_sigma} for all $k \in \N{}_+$, in which case the Lipschitz constant estimates are monotonically increasing.  However, in Algorithm~\ref{alg.dynamic}, we allow these estimates to decrease since this might yield better results in practice.
  \item If the Lipschitz constants $L$ and $\sigma$ for $g := \nabla f$ and $H := \nabla^2 f$, respectively, are known, then one could simply set $L_k = L$ and $\sigma_k = \sigma$ during each iteration; in this case, the \textbf{loop} would not actually be needed.  Although this would simplify the presentation, it would generally result in more iterations being required to obtain (approximate) first- and/or second-order stationarity.  However, if the cost of evaluating $f$ is substantial, then static parameters might work well.   
  \item Each time through the \textbf{loop}, condition~\eqref{mod.dec.alpha} or~\eqref{mod.dec.beta} is tested, but not both, since this would require an extra evaluation of $f$. If the cost of evaluating the objective function is not a concern, then one could choose between the two steps based on \emph{actual} objective function decrease rather than model decrease.
\eitemize

%************
% Subsection
%************
\subsection{Complexity Analysis}\label{sec.deterministic.complexity}

We have the following complexity result for Algorithm~\ref{alg.dynamic}.  We claim that a similar result could be stated and proved for Algorithm~\ref{alg.2step} as well, where one employs the inequality \eqref{eq.twostep_reduction} in place of \eqref{dec.dynamic} in the proof below.  However, for brevity and since we believe it might often be the better method in practice, we focus on the following result for our dynamic method, Algorithm~\ref{alg.dynamic}.

\btheorem\label{thm:complexity}
  Consider any scalars $\epsilon_g \in (0,\bar\epsilon_g]$ and $\epsilon_H \in (0,\bar\epsilon_H]$ for some constants $(\bar\epsilon_g,\bar\epsilon_H) \in (0,\infty) \times (0,\infty)$.  With respect to Algorithm~\ref{alg.dynamic}, the cardinality of the index set
  \bequationNN
    \Gcal(\epsilon_g) := \{k \in \N{}_+ : \|g(x_k)\|_2 > \epsilon_g\}
  \eequationNN
  is at most $\Ocal(\epsilon_g^{-2})$.  In addition, the cardinality of the index set
  \bequationNN
    \Hcal(\epsilon_H) := \{k \in \N{}_+ : |(\lambda_k)_-| > \epsilon_H\}
  \eequationNN
  is at most $\Ocal(\epsilon_H^{-3})$.  Hence, the number of iterations and derivative (i.e., gradient and Hessian) evaluations required until iteration $k \in \N{}_+$ is reached with
  \bequationNN
    \|g(x_k)\|_2 \leq \epsilon_g\ \ \text{and}\ \ |(\lambda_k)_-| \leq \epsilon_H
  \eequationNN
  is at most $\Ocal(\max\{\epsilon_g^{-2},\epsilon_H^{-3}\})$.  Moreover, if Step~\ref{step.set_L_sigma} is modified such that
  \bequation\label{eq.L_sigma_lower}
    (L_{k+1},\sigma_{k+1}) \in [L_{\min},L_k] \times [\sigma_{\min},\sigma_k],
  \eequation
  for some $(L_{\min},\sigma_{\min}) \in (0,\infty) \times (0,\infty)$ with $L_1 \geq L_{\min}$ and $\sigma_1 \geq \sigma_{\min}$, then the number of iterations in the \textbf{loop}  for all $k \in \N{}_+$ is uniformly bounded, meaning that the complexity bound above also holds for the number of function evaluations.
\etheorem
\bproof
  As in the proof of Theorem~\ref{th.deterministic.adaptive}, one has that, for all $k \in \N{}_+$, the inequality~\eqref{dec.dynamic} holds, which we restate here as
  \bequation\label{dec.dynamic_again}
    f(x_k) - f(x_{k+1}) \geq \max\left\{\frac{\delta^2}{2L_k} \|g(x_k)\|_2^2, \frac{2\gamma^3}{3\sigma_k^2} |(\lambda_k)_-|^3 \right\}.
  \eequation
  From this inequality and the bounds in \eqref{eq.L_sigma_upper}, it follows that
  \bequationNN
    \baligned
      k \in \Gcal(\epsilon_g) &\implies f(x_k) - f(x_{k+1}) \geq \frac{\delta^2}{2L_{\max}} \epsilon_g^2 \\
      \text{while}\ \ 
      k \in \Hcal(\epsilon_H) &\implies f(x_k) - f(x_{k+1}) \geq \frac{2\gamma^3}{3\sigma_{\max}^2} \epsilon_H^3.
    \ealigned
  \eequationNN
  Since $f$ is bounded below by $f_{\inf}$ and \eqref{dec.dynamic_again} ensures that $\{f(x_k)\}$ is monotonically decreasing, the inequalities above imply that $\Gcal(\epsilon_g)$ and $\Hcal(\epsilon_H)$ are both finite.  In addition, by summing the reductions achieved in $f$ over the iterations in the above index sets, we obtain
  \bequationNN
    \baligned
      f(x_1) - f_{\inf} &\geq \sum_{k\in\Gcal(\epsilon_g)} f(x_k) - f(x_{k+1}) \geq |\Gcal(\epsilon_g)| \frac{\delta^2}{2L_{\max}} \epsilon_g^2 \\
      \text{and}\ \ 
      f(x_1) - f_{\inf} &\geq \sum_{k\in\Hcal(\epsilon_H)} f(x_k) - f(x_{k+1}) \geq |\Hcal(\epsilon_H)| \frac{2\gamma^3}{3\sigma_{\max}^2} \epsilon_H^3.
    \ealigned
  \eequationNN
  Rearranging, one obtains that
  \bequationNN
    |\Gcal(\epsilon_g)| \leq \(\frac{2L_{\max}(f(x_1) - f_{\inf})}{\delta^2}\)\epsilon_g^{-2}\ \ \text{and}\ \ |\Hcal(\epsilon_H)| \leq \(\frac{3\sigma_{\max}^2(f(x_1) - f_{\inf})}{2\gamma^3}\)\epsilon_H^{-3},
  \eequationNN
  as desired.  Finally, the last conclusion follows by the fact that, with~\eqref{eq.L_sigma_lower}, the number of iterations within the \textbf{loop} for any $k \in \N{}_+$ is bounded above by
  \bequationNN
    \(1 + \left\lfloor\frac{1}{\log(\rho)}\log\(\frac{L}{L_{\min}}\)\right\rfloor\) + \(1 + \left\lfloor\frac{1}{\log(\rho)}\log\(\frac{\sigma}{\sigma_{\min}}\)\right\rfloor\),
  \eequationNN
  which is the maximum number of updates to $L_k$ and/or $\sigma_k$ to bring values within $[L_{\min},L]$ and $[\sigma_{\min},\sigma]$ up to satisfy $L_k \geq L$ and $\sigma_k \geq \sigma$. \qed
\eproof

%************
% Subsection
%************
\subsection{Behavior Near Strict Saddle Points}\label{sec.deterministic.discussion}

As previously mentioned (recall \S\ref{sec.prior}), much recent attention has been directed toward the behavior of nonconvex optimization algorithms in the neighborhood of saddle points.  In particular, in order to prove guarantees about avoiding saddle points, an assumption is often made about all saddle points being \emph{strict} (or \emph{ridable}) or even \emph{non-degenerate}.  Strict saddle points are those at which the Hessian has at least one negative eigenvalue; intuitively, these are saddle points from which an algorithm should usually be expected to avoid.  Nondegenerate saddle points are ones at which the Hessian has no zero eigenvalues, a much stronger assumption.

One can show that in certain problems of interest, all saddle points are strict.  This is interesting.  However, much of the work that has discussed the behavior of nonconvex algorithms in the neighborhood of strict saddle points have focused on standard types of descent methods, about which one can only prove high probability results \cite{LeePanPilSimJorRec17}.  By contrast, when one considers explicitly computing directions of negative curvature, this is all less of an issue.  After all, notice that our convergence and complexity analyses for our framework did not require careful consideration of the nature of any potential saddle points.

In any case, for ease of comparison to other recent analyses, let us discuss the convergence/complexity properties for our method in the context of a strict saddle point assumption.  Suppose that the set of maximizers and saddle points of~$f$, say $\{\bar x_i\}_{i\in\cal{I}}$ for some index set~$\cal{I}$, which must necessarily have $g(\xbar_i) = 0$ for all $i \in \Ical$, also has $\bar\lambda_i < 0$ for all $i \in \Ical$, where $\{\bar\lambda_i\}_{i\in\Ical}$ are the leftmost eigenvalues of $\{H(\bar x_i)\}_{i\in\Ical}$ and are uniformly bounded away from zero. In this setting, we may draw the following conclusions from Theorems~\ref{th.deterministic.adaptive} and~\ref{thm:complexity}.
\bitemize
\item Any limit point of the sequence $\{x_k\}$ computed by Algorithm~\ref{alg.dynamic} is a minimizer of $f$.  This follows since Theorem~\ref{th.deterministic.adaptive} ensures that for any limit point, say $\xbar$, the gradient must be zero and the leftmost eigenvalue must be nonnegative.  It follows from these facts and the assumptions on the maximizers and saddle points $\{\xbar_i\}$ that $\xbar$ must be a minimizer, as claimed. 
\item  From the discussion in the previous bullet and Theorem~\ref{thm:complexity} we know, in fact, that the iterates of Algorithm~\ref{alg.dynamic} must eventually enter a region consisting only of minimizers (i.e., one that does not contain any maximizers or saddle points) in a number of iterations that is polynomial in $\epsilon \in (0,\infty)$, where the negative value $-\epsilon$ is greater than the largest of the leftmost eigenvalues of $\{H(\xbar_i)\}_{i\in\Ical}$. This complexity result holds without assuming that the minimizers, maximizers, and saddle points are all nondegenerate (i.e., the Hessian matrix at these points are nonsingular), which, e.g., should be contrasted with the analysis in~\cite[see Assumption~3]{PateMokhRibe17}.  The primary reason for our stronger convergence/complexity properties is that our method incorporates negative curvature directions when they exist, as opposed to only descent directions.
\eitemize

%\textbf{NEED TO DO.}\footnote{dpr: if we use the results when sigma/L are only allowed to increase, then our results are better in all cases. Can also talk about lack of fast local convergence.}

%************
% Subsection
%************
\subsection{Step Computation Techniques}\label{sec.deterministic.step_computation}

There is flexibility in the ways in which the steps $d_k$ and~$s_k$ (or $\shat_k$) are computed in order to satisfy the desired conditions (in \eqref{def.d} and \eqref{def.s}/\eqref{def.s.adapt}).  For example, $s_k$ might be the steepest descent direction $-g(x_k)$, for which \eqref{def.s} holds with $\delta = \zeta = \eta = 1$.  Another option, with a symmetric positive definite $B_k^{-1} \in \R{n \times n}$ that has $\kappa(B_k^{-1}) \leq \delta^{-1}$ with a spectrum falling in $[\zeta,\eta]$, is to compute $s_k$ as the modified Newton direction $-B_k^{-1}g(x_k)$.  A particularly attractive option for certain applications (when Hessian-vector products are readily computed) is to compute~$s_k$ via a Newton-CG routine \cite{NoceWrig06} with safeguards that terminate the iteration before a CG iterate is computed that violates~\eqref{def.s} for prescribed $(\delta,\zeta,\eta)$.

There are multiple ways to compute the negative curvature direction.  In theory, the most straightforward approach is to set $d_k = \pm v_k$, where $v_k$ is a leftmost eigenvector of $H(x_k)$.  With this choice, it follows that $d_k^T H(x_k) d_k = \lambda_k \|d_k\|^2$, meaning that~\eqref{def.d.curv} is satisfied, and one can scale $d_k$ to ensure that~\eqref{def.d.descent} and \eqref{def.d.norm} hold.  A second approach for large-scale settings (i.e., when $n$ is large) is to compute $d_k$ via matrix-free Lanczos iterations~\cite{lanczos1950iteration}.  Such an approach can produce a direction $d_k$ satisfying~\eqref{def.d.curv}, which can then be scaled to yield~\eqref{def.d.descent} and \eqref{def.d.norm}.

%************
% Subsection
%************
\subsection{Numerical Results}\label{sec.deterministic.numerical}

In this section, we demonstrate that there can be practical benefits of following directions of negative curvature if one follows the dynamic approach of Algorithm~\ref{alg.dynamic}.  To do this, we implemented software in {\sc Matlab} that, for all $k \in \N{}_+$, has the option to compute $s_k$ via several options (see \S\ref{sec.deterministic.sd} and \S\ref{sec.deterministic.modified}) and $d_k = \pm v_k$ (recall \S\ref{sec.deterministic.step_computation}).  Our test problems include a subset of the {\sc CUTEst} collection~\cite{gould2015cutest}.  Specifically, we selected  all of the unconstrained problems with $n \leq 500$ and second derivatives explicitly available.  This left us with a test set of $97$ problems.
%unconstrained optimization problems.

%  Second, we ran Algorithm~\ref{alg.dynamic} on this subset for $100$ iterations and only kept those for which at least one direction of negative curvature was used.  This process left us with the $26$ test problems listed under the column ``Problem'' in Table~\ref{tab.deterministic}.

Using this test set, we considered two variants of Algorithm~\ref{alg.dynamic}: $(i)$ a version in which the \textbf{if} condition in Step~7 is always presumed to test true, i.e., the descent step~$s_k$ is chosen for all $k \in \N{}_+$ (which we refer to as Algorithm~\ref{alg.dynamic}($s_k$)), and $(ii)$ a version that, as in our formal statement of the algorithm, chooses between descent and negative curvature steps by comparing model reduction values for each $k \in \N{}_+$ (which we refer to as Algorithm~\ref{alg.dynamic}($s_k,d_k$)). In our experiments, we used $L_1 \gets 1$ and $\sigma_1 \gets 1$, updating them and setting subsequent values in their respective sequences using the following strategy.  For increasing one of these values in Step~\ref{step.L-increase} or Step~\ref{step.sigma-increase}, we respectively set the quantities
\begin{align*}
\hat L_k 
&\gets L_k + \frac{2\big(f(x_k+\alpha_k s_k)-f(x_k) + m_{s,k}(\alpha_k)\big)}{\alpha_k^2 \|s_k\|^2} \ \ \text{or}  \\
\hat \sigma_k 
&\gets \sigma_k + \frac{6\big(f(x_k+\beta_k d_k) - f(x_k) + m_{d,k}(\beta_k)\big)}{\beta_k^3 \|d_k\|^3},
\end{align*}
then, with $\rho \gets 2$, use the update
\begin{equation}\label{update.L.and.sigma}
\begin{aligned}
L_k 
&\gets  \max\{\rho L_k, \min\{10^{3}L_k,\hat L_k\}\} 
\ \ \text{in Step~\ref{step.L-increase} of Algorithm~\ref{alg.dynamic} or} \\
%&\qquad \text{(Step~\ref{step.L-increase} of Algorithm~\ref{alg.dynamic})} \\
\sigma_k 
&\gets \max\{\rho\sigma_k,\min\{10^3\sigma_k,\hat\sigma_k\}\}
 \ \ \text{in Step~\ref{step.sigma-increase} of Algorithm~\ref{alg.dynamic}.}
%&\qquad \text{(Step~\ref{step.sigma-increase} of Algorithm~\ref{alg.dynamic})}
\end{aligned}
\end{equation}
The quantity $\hat L_k$ is defined so that \eqref{mod.dec.alpha} holds at equality with $\alpha = \alpha_k$ and $L_k$ replaced by $\hat L_k$ in the definition of $m_{s,k}(\alpha_k)$, i.e., $\hat L_k$ is the value that makes the model decrease agree with the exact function decrease at $x_k + \alpha_k s_k$.  The procedure we use to set $\hat \sigma_k$ is analogous.  We remark that the updates in~\eqref{update.L.and.sigma} ensure that $L_k \in [\rho L_k,10^3L_k]$ and $\sigma_k \in [\rho\sigma_k,10^3\sigma_k]$, which we claim maintains the convergence and complexity guarantees of Algorithm~\ref{alg.dynamic}.  Moreover, when Step~\ref{step.L-increase} (respectively, Step~\ref{step.sigma-increase}) is reached, then it must be the case that $\hat L_k > L_k$ (respectively, $\hat\sigma_k > \sigma_k$).  On the other hand, in Step~\ref{step.set_L_sigma} we use the following updates:
\begin{equation*}%\label{update.L.and.sigma.decrease}
\begin{aligned}
L_{k+1} 
&\gets
\max\{10^{-3},10^{-3}L_k, \hat L_k\}
\ \ \text{and} \ \ \sigma_{k+1} \gets \sigma_k 
\ \ \text{when}\ \ x_{k+1} \gets x_k + \alpha_k s_k; \\
\sigma_{k+1} 
&\gets
\max\{10^{-3},10^{-3}\sigma_k, \hat \sigma_k\}
\ \ \text{and} \ \ L_{k+1} \gets L_k
\ \ \text{when}\ \ x_{k+1} \gets x_k + \beta_k d_k.
\end{aligned}
\end{equation*}
%In Step~\ref{step.set_L_sigma}, we set $L_{k+1} \gets 0.5 L_k$ if $s_k$ has been taken or set $\sigma_{k+1} \gets 0.5\sigma_k$ if $d_k$ has been taken.  
Over the next two sections, we discuss numerical results when two different options for computing the descent direction $s_k$ are used.

\subsubsection{Choosing $s_k$ as the direction of steepest descent} \label{sec.deterministic.sd}

The tests in this section use the steepest descent direction, i.e., $s_k = -g(x_k)$ for all $k\in\N{}_+$. Although this is a simple choice, it gives a starting point for understanding the potential benefits of using directions of negative curvature. 

\begin{figure}[ht]
\centering
\begin{subfigure}[b]{0.47\textwidth}
                \centering
                \includegraphics[scale=0.36,angle=-0,trim={2.8cm 7cm 2.5cm 7cm},clip]{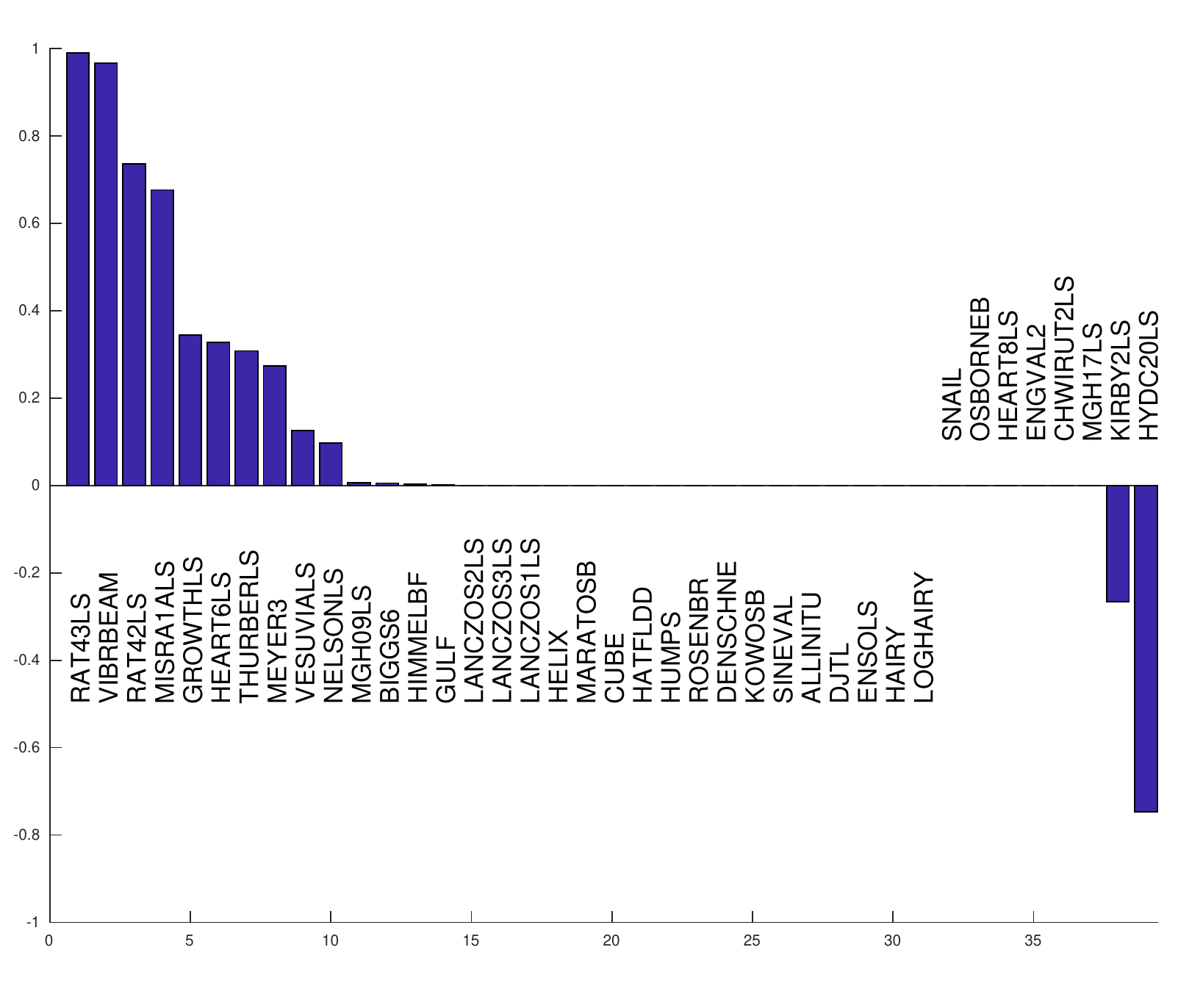}
                \caption{Plot associated with the quantity~\eqref{measure.f_diff}.} %for the choice $s_k = -g_k$ 
                %for all test problems that used at least one negative curvature direction.}
                \label{fig:sd.term.f_diff}
\end{subfigure} \\
\begin{subfigure}[b]{0.47\textwidth}
                \centering
                \includegraphics[scale=0.36,angle=-0,trim={2.8cm 7cm 2.5cm 7cm},clip]{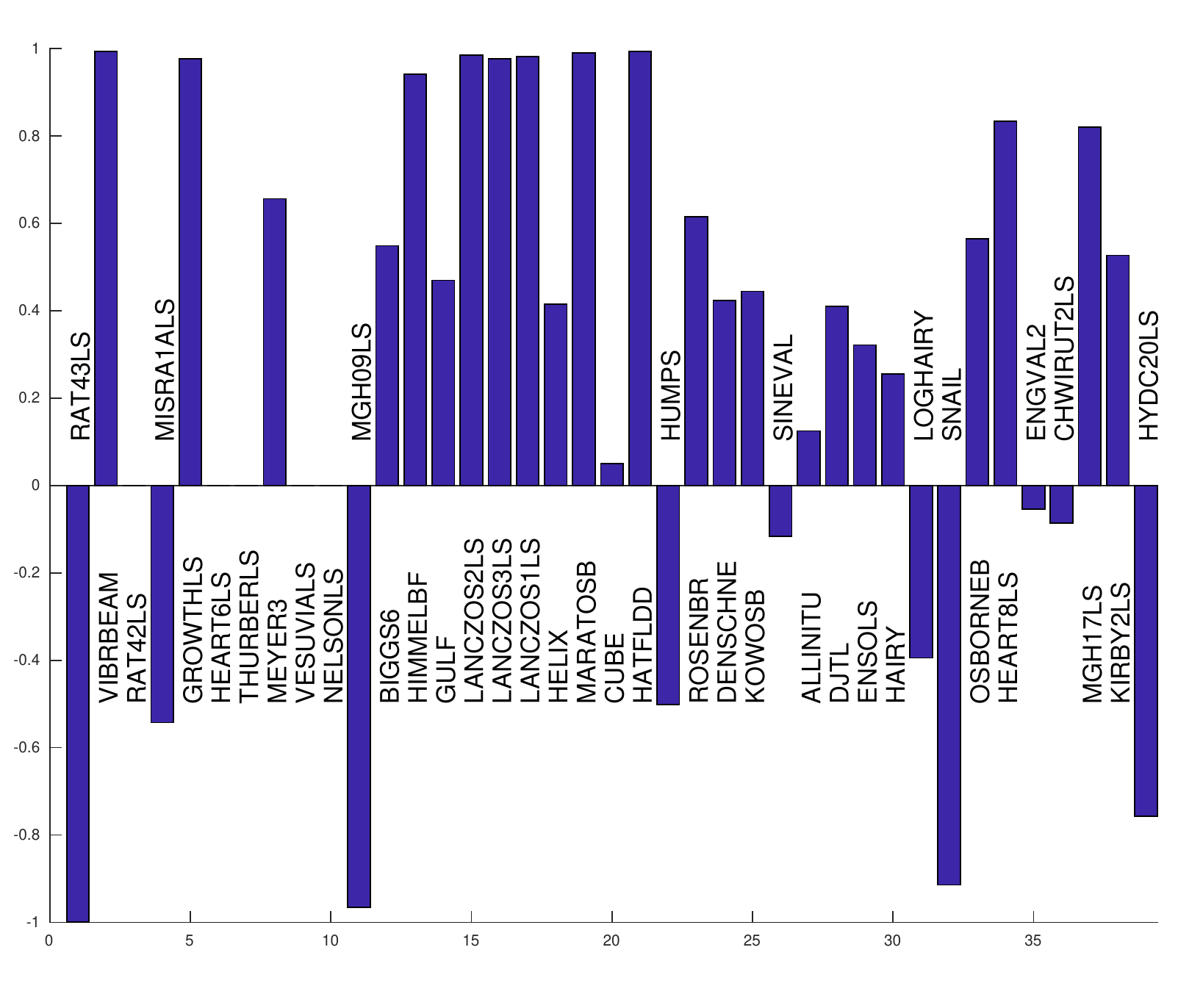}
                \caption{Plot associated with the quantity~\eqref{measure.iterates}.} %for the choice $s_k = -g_k$ 
                %for all test problems that used at least one negative curvature direction.}
                \label{fig:sd.term.iterates}
\end{subfigure} 
\phantom{--}
\begin{subfigure}[b]{0.47\textwidth}
                \centering
                \includegraphics[scale=0.36,angle=-0,trim={2.8cm 7cm 2.5cm 7cm},clip]{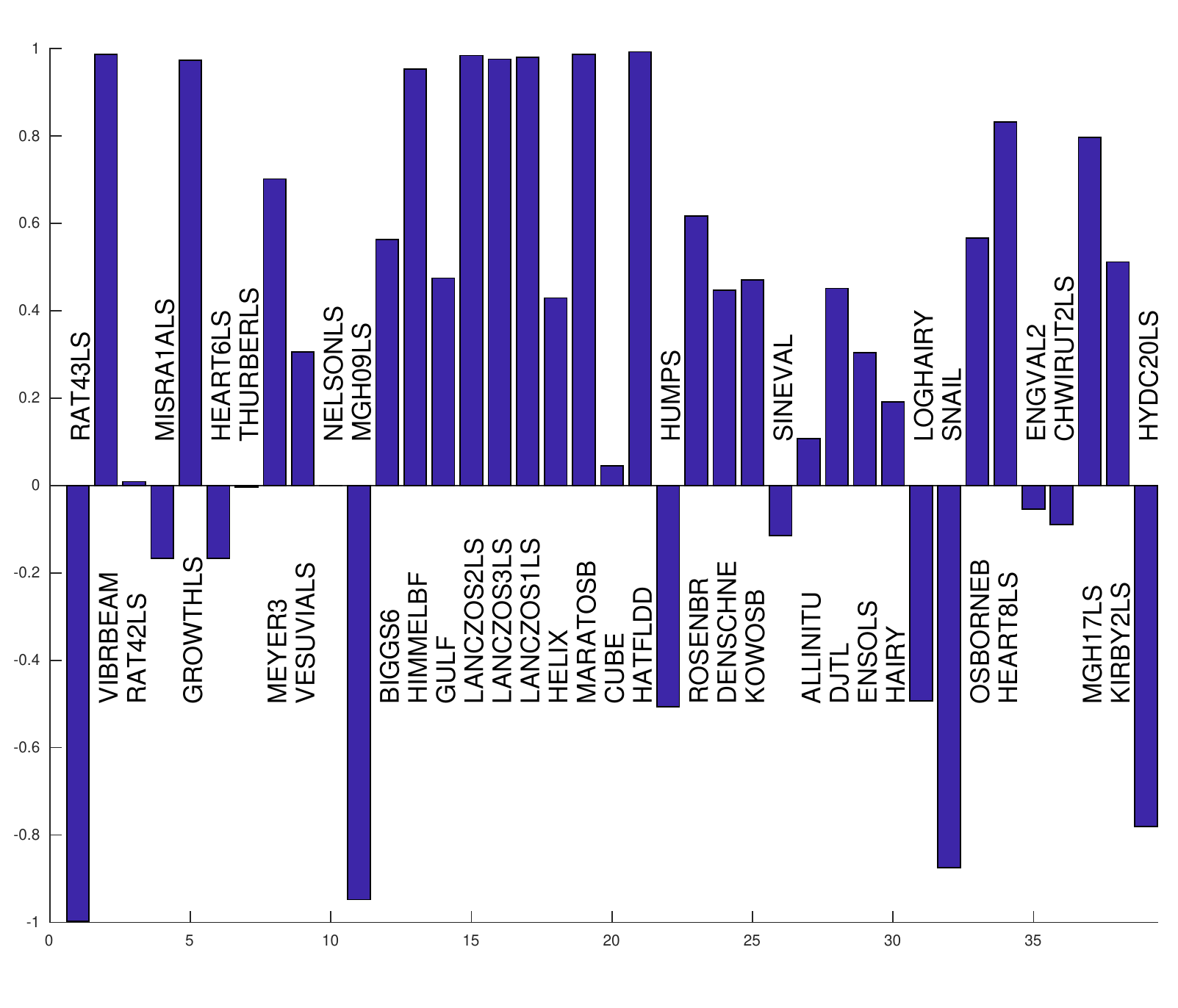}
                \caption{Plot associated with the quantity~\eqref{measure.fevals}.} %for the choice $s_k = -g_k$ 
                %for all test problems that used at least one negative curvature direction.}
                \label{fig:sd.term.fevals}
\end{subfigure}
\caption{Plots for the choice $s_k \equiv -g(x_k)$ for all $k \in \N{}_+$.
% with termination condition~\eqref{term.2} used for Algorithm~\ref{alg.dynamic}($s_k$) and Algorithm~\ref{alg.dynamic}($s_k,d_k$), respectively. 
Only problems for which at least one negative curvature direction is used 
%and for which the value in~\eqref{measure.f_diff} is larger than $10^{-5}$ in absolute value 
are presented.  The problems are ordered based on the values in plot (a).}\label{fig:sd.term.all}
%these bar plots display information related to the final objective value (Figure~\ref{fig:sd.term.f_diff}), the number of required iterations (Figure~\ref{fig:sd.term.iterates}), and the number of objective function evaluations (Figure~\ref{fig:sd.term.fevals}), for all test problems for which at least one negative curvature direction was computed.}\label{fig:sd.term.all}
\end{figure}

We ran Algorithm~\ref{alg.dynamic}($s_k$) and Algorithm~\ref{alg.dynamic}($s_k,d_k$) on the previously described {\sc CUTEst} problems with commonly used stopping conditions.  Specifically, an algorithm terminates with an approximate second-order stationary solution if
%\begin{equation} \label{term.1}
%\|g(x_k)\| \leq 10^{-5} \max\{1, \|g(x_0)\|\},
%\end{equation}
%and Algorithm~\ref{alg.dynamic}($s_k$,$d_k$) terminated with an approximate second-order solution if
\begin{equation} \label{term.2}
\|g(x_k)\| \leq 10^{-5} \max\{1, \|g(x_1)\|\}
\ \ \text{and} \ \
|(\lambda_k)_-| \leq 10^{-5} \max\{1, |(\lambda_1)_-|\}.
\end{equation}
We also terminate an algorithm if an iteration limit of $10,\!000$ is reached or if a trial step smaller than $10^{-16}$ is computed. The results can be found in Figure~\ref{fig:sd.term.all}.

%============================================
% Note: do not erase
% This data came from directory analyze_results_output_27.
 %============================================

In Figure~\ref{fig:sd.term.f_diff}, letting $f_{\text{final}}(s_k)$ and $f_{\text{final}}(s_k,d_k)$ be the final computed objective values for Algorithm~\ref{alg.dynamic}($s_k$) and Algorithm~\ref{alg.dynamic}($s_k,d_k$), respectively, we plot
\begin{equation}\label{measure.f_diff}
\frac{f_{\text{final}}(s_k) - f_{\text{final}}(s_k,d_k)}{\max\{|f_{\text{final}}(s_k)|,|f_{\text{final}}(s_k,d_k)|,1\}}
\in [-2,2]
\end{equation}
for each problem that Algorithm~\ref{alg.dynamic}($s_k,d_k$) used at least one negative curvature direction.
% and for which~\eqref{measure.f_diff} is larger than $10^{-5}$ in absolute value.  
In this manner, an upward pointing bar implies that  Algorithm~\ref{alg.dynamic}($s_k,d_k$) terminated with a lower value of the objective function, while its magnitude represents how much better was this value. We can observe from Figure~\ref{fig:sd.term.f_diff} that among the $39$ problems,
% that a direction of negative curvature was used and for which  
Algorithm~\ref{alg.dynamic}($s_k$) terminated with a significantly lower objective value compared to Algorithm~\ref{alg.dynamic}($s_k,d_k$) only $2$ times.  (While the value in \eqref{measure.f_diff} falls in the interval $[-2,2]$, the bars in the figure all fall within $[-1,1]$ since, for each problem, the final function value reached by each algorithm had the same sign.)

We are also interested in the number of iterations and objective function evaluations needed to obtain these final objective function values.  The relative performances of the algorithms with respect to these measures are shown in Figure~\ref{fig:sd.term.iterates} and Figure~\ref{fig:sd.term.fevals}.  In particular, letting $\#its({s_k})$ be the total number of iterations required by Algorithm~\ref{alg.dynamic}($s_k$) and $\#its({s_k,d_k})$ be the total number of iterations required by Algorithm~\ref{alg.dynamic}($s_k,d_k$), we plot in Figure~\ref{fig:sd.term.iterates} the values
\begin{equation}\label{measure.iterates}
\frac{\#its({s_k}) - \#its({s_k,d_k)}}{\max\{\#its({s_k}),\#its({s_k,d_k}),1\}}
\in [-1,1],
\end{equation}
and, letting $\#fevals({s_k})$ be the total number of objective function evaluations for Algorithm~\ref{alg.dynamic}($s_k$) and $\#fevals({s_k,d_k})$ be the total number of objective function evaluations for Algorithm~\ref{alg.dynamic}($s_k,d_k$), we plot in Figure~\ref{fig:sd.term.fevals} the values 
\begin{equation}\label{measure.fevals}
\frac{\# fevals({s_k}) - \# fevals({s_k,d_k})}{\max\{\#fevals({s_k}),\#fevals({s_k,d_k}),1\}}
\in [-1,1].
\end{equation}
These two plots show that Algorithm~\ref{alg.dynamic}($s_k,d_k$) tends to perform better than Algorithm~\ref{alg.dynamic}($s_k$) in terms of both the number of required iterations and function evaluations. Overall, we find these results interesting since Algorithm~\ref{alg.dynamic}($s_k,d_k$) tends to find lower values of the objective function (see Figure~\ref{fig:sd.term.f_diff}) while typically requiring fewer iterations (see Figure~\ref{fig:sd.term.iterates}) and function evaluations (see Figure~\ref{fig:sd.term.fevals}).

\subsubsection{Choosing $s_k$ using a modified-Newton strategy} \label{sec.deterministic.modified}

In this subsection, we show the results of tests similar to those in
\S\ref{sec.deterministic.sd}
% for the previous subsection 
except that now we compute the descent direction $s_k$ by a modified-Newton approach.  Specifically, we compute $s_k$ as the unique vector satisfying $B_k s_k = -g(x_k)$ with 
\begin{equation} \label{def:Bk}
B_k = H(x_k) + \delta_k I,
\end{equation}
where $I$ is the identify matrix and $\delta_k$ is the smallest nonnegative real number such that $B_k$ is positive definite with a condition number less than or equal to $10^8$.
%\footnote{FEC: Do you think it's worth making a comment here that since we compute $\delta_k$ this way, computing the negative curvature direction is not much extra work, so the two algorithms have similar per-iteration costs here?}

\begin{center}
\begin{figure}[t!]
\centering
\begin{subfigure}[b]{0.47\textwidth}
                \centering
                \includegraphics[scale=0.36,angle=-0,trim={2.8cm 7cm 2.5cm 7cm},clip]{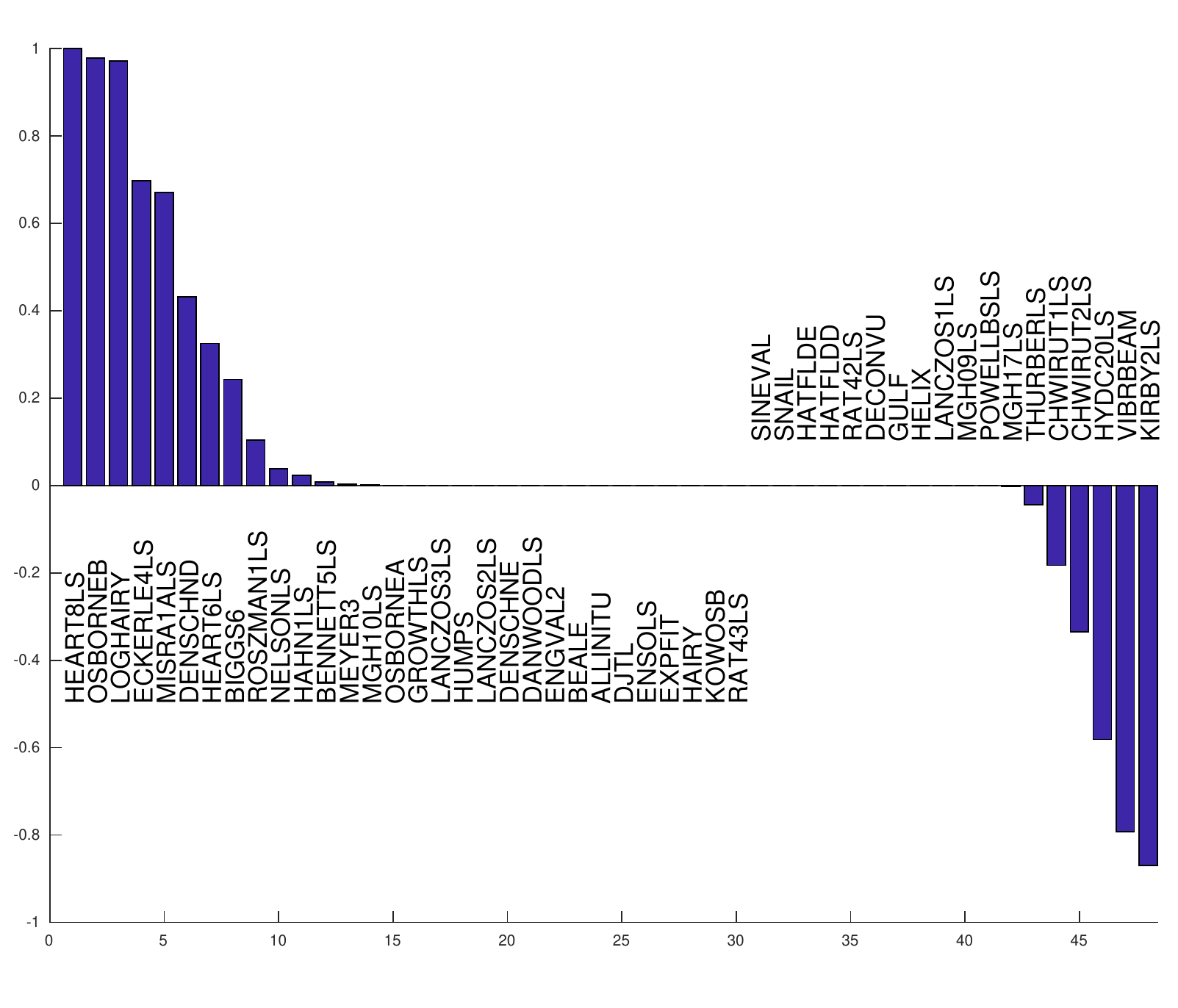}
                \caption{Plot associated with the quantity~\eqref{measure.f_diff}.}%for the choice $s_k = -g_k$ 
                %for all test problems that used at least one negative curvature direction.}
                \label{fig:shifted.term.f_diff}
\end{subfigure} \\
\begin{subfigure}[b]{0.47\textwidth}
                \centering
                \includegraphics[scale=0.36,angle=-0,trim={2.8cm 7cm 2.5cm 7cm},clip]{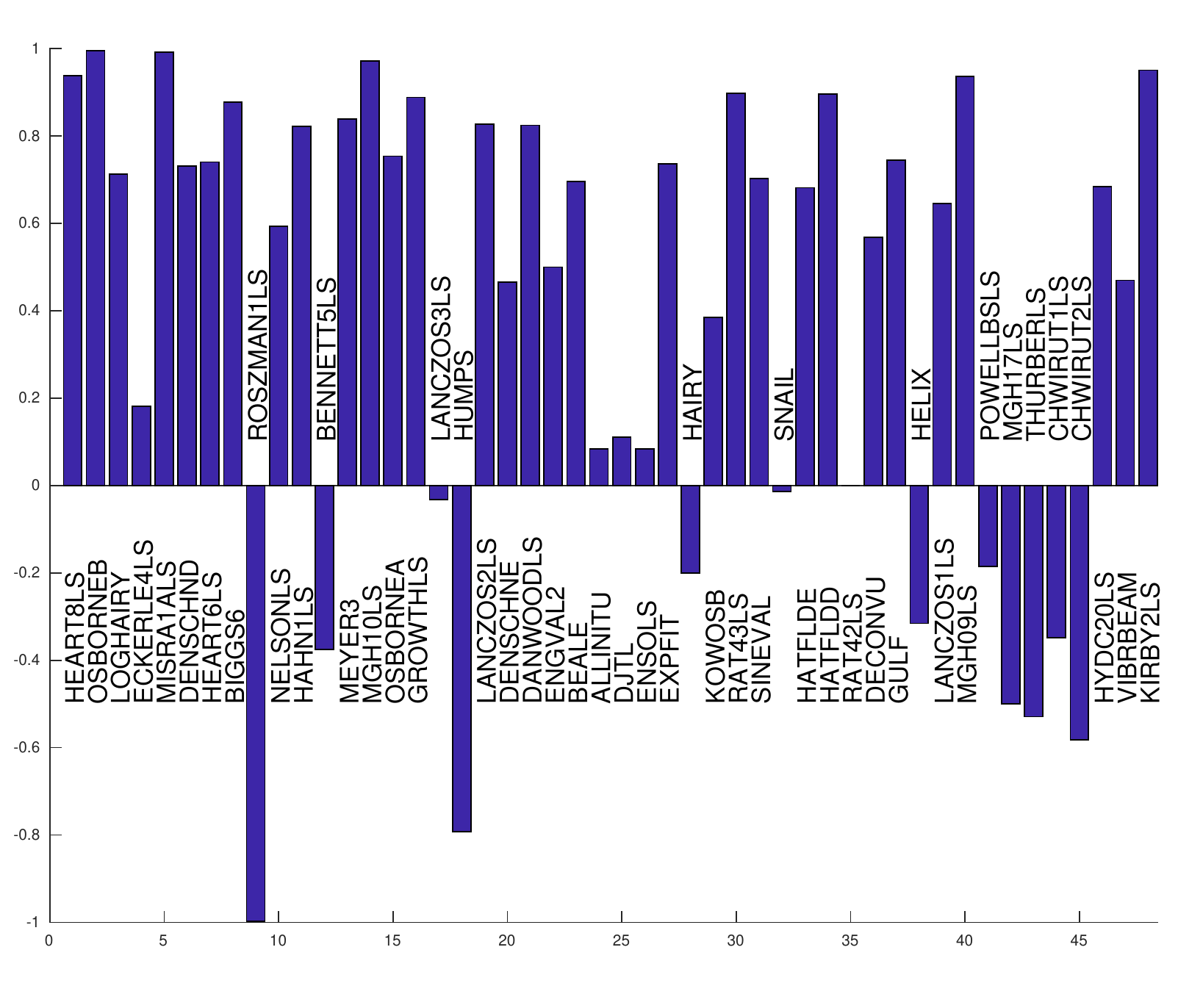}
                \caption{Plot associated with the quantity~\eqref{measure.iterates}.} %for the choice $s_k = -g_k$ 
                %for all test problems that used at least one negative curvature direction.}
                \label{fig:shifted.term.iterates}
\end{subfigure}
\phantom{--}
\begin{subfigure}[b]{0.47\textwidth}
                \centering
                \includegraphics[scale=0.36,angle=-0,trim={2.8cm 7cm 2.5cm 7cm},clip]{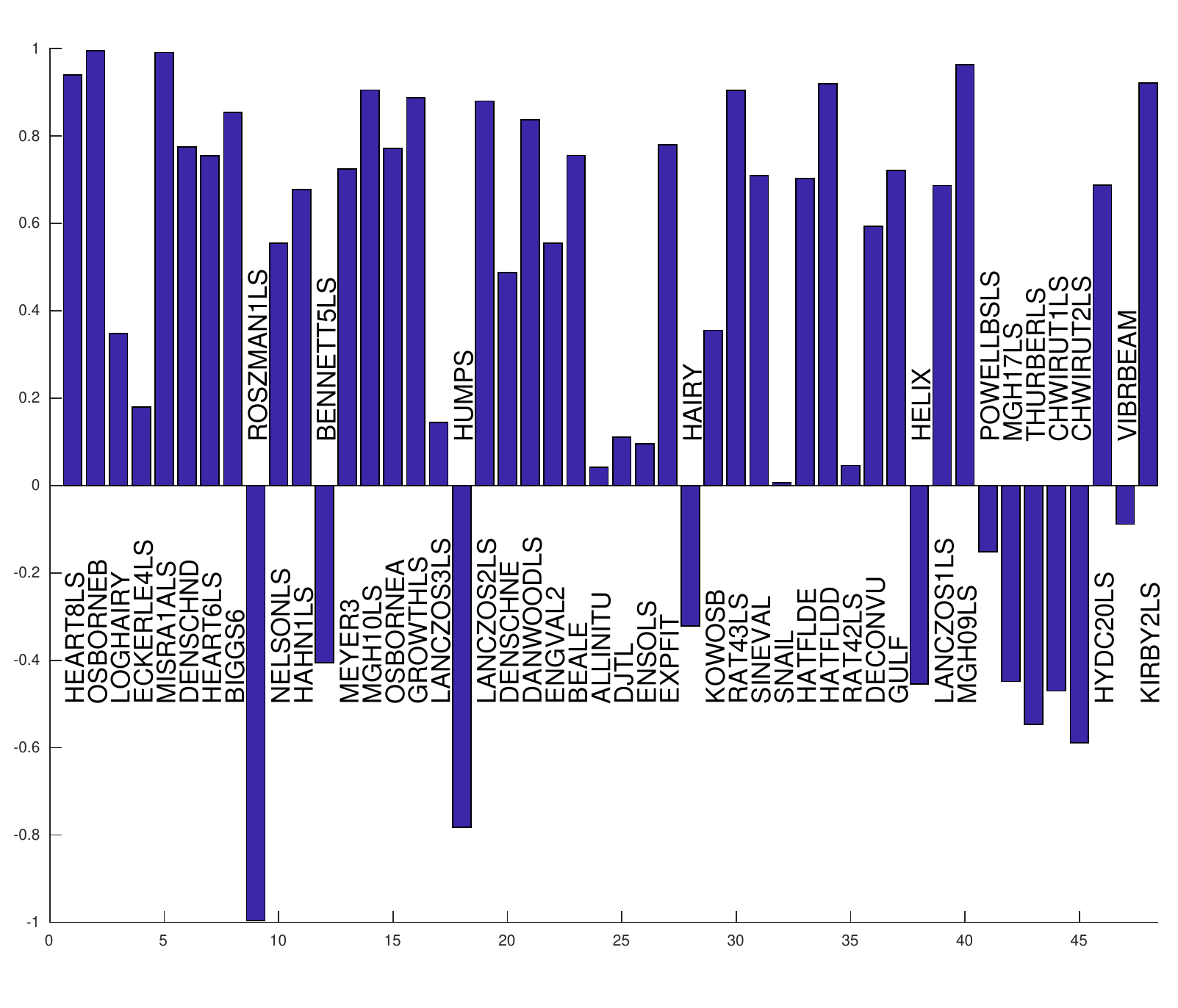}
                \caption{Plot associated with the quantity~\eqref{measure.fevals}.} %for the choice $s_k = -g_k$ 
                %for all test problems that used at least one negative curvature direction.}
                \label{fig:shifted.term.fevals}
\end{subfigure}
\caption{Plots for when $s_k$ is chosen to satisfy $B_k s_k = -g(x_k)$ with $B_k$ defined via~\eqref{def:Bk}.
% and when termination conditions~\eqref{term.1} and~\eqref{term.2} are used for Algorithm~\ref{alg.dynamic}($s_k$) and Algorithm~\ref{alg.dynamic}($s_k,d_k$), respectively. 
Only problems for which at least one negative curvature direction is used 
%and for which the value in~\eqref{measure.f_diff} is larger than $10^{-5}$ in absolute value 
are presented.  The problems are ordered based on the values in plot (a).}\label{fig:shifted.term.all}
%these bar plots display information related to the final objective value (Figure~\ref{fig:sd.term.f_diff}), the number of required iterations (Figure~\ref{fig:sd.term.iterates}), and the number of objective function evaluations (Figure~\ref{fig:sd.term.fevals}), for all test problems for which at least one negative curvature direction was computed.}\label{fig:sd.term.all}
\end{figure}
\end{center}

The results from solving our {\sc CUTEst} test set with this choice of $B_k$ using the stopping condition~\eqref{term.2}
% for Algorithm~\ref{alg.dynamic}($s_k$) and the stopping condition~\eqref{term.2} for Algorithm~\ref{alg.dynamic}($s_k,d_k$) 
are presented in Figure~\ref{fig:shifted.term.all}. We can observe that, overall, Algorithm~\ref{alg.dynamic}($s_k,d_k$) outperforms Algorithm~\ref{alg.dynamic}($s_k$) in this experiment.  Neither algorithm consistently outforms the other in terms of their final objective values, but Algorithm~\ref{alg.dynamic}($s_k,d_k$) typically requires fewer iterations (see Figure~\ref{fig:shifted.term.iterates}) and objective function evaluations (see Figure~\ref{fig:shifted.term.fevals}).

%Finally, we ran the test described in the previous section associated with the quantity~\eqref{term.fval}. Namely, we ran both algorithms and terminated when the objective function reached the target value~\eqref{term.fval}. 
%%, i.e., we ran both algorithms until they reached $95\%$ of the reduction in the objective function obtained by the worst of the two methods.  
%The results, which may be found in Figure~\ref{fig:shifted.fval.iterates_and_fevals}, illustrate that directions of negative curvature have no advantage in terms of the number of required iterations and seems to be disadvantageous in terms of the number of required objective function evaluations, for this experimental setup.

%at least in the current setting where $B_k s_k = -g_k$ and $B_k$ is computed via the modified Newton strategy in~\eqref{def:Bk}.

%============================================
% Note: do not erase
% This data came from directory analyze_results_output_27.
 %============================================

\section{Stochastic Optimization}\label{sec.stochastic}

Let us now consider the problem to minimize a twice continuously differentiable and bounded below (by $f_{\inf} \in \R{}$) objective function $f : \R{n} \to \R{}$ defined by the expectation, in terms of the distribution of a random variable $\xi$ with domain $\Xi$, of a stochastic function $F : \R{n} \times \Xi \to \R{}$, namely,
\bequation\label{prob.expectation}
  \min_{x \in \R{n}}\ f(x),\ \ \text{where}\ \ f(x) := \E_\xi[F(x,\xi)].
\eequation
In this context, we expect that, at an iterate $x_k$, one can only compute stochastic gradient and Hessian estimates.  We do not claim that we are able to prove convergence guarantees to second-order stationarity as in the deterministic case.  That said, we are able to present a two-step method with convergence guarantees to first-order stationarity whose structure motivates a dynamic method that we show can offer beneficial practical performance by exploring negative curvature.

%************
% Subsection
%************
\subsection{Two-Step Method: Stochastic Gradient/Newton with ``Curvature Noise''}\label{sec.stochastic.two-step}

At an iterate $x_k$, let $\xi_k$ be a random variable representing a seed for generating a vector $s_k \in \R{n}$.  For example, if $f$ is the expected function value over inputs from a dataset, then~$\xi_k$ might represent sets of points randomly drawn from the dataset.  With $\E_{\xi_k}[\cdot]$ denoting expectation taken with respect to the distribution of $\xi_k$ given the current iterate $x_k$, we require the vector $s_k$ to satisfy
\bsubequations\label{def.stochastic.s}
  \begin{align}
    -\nabla f(x_k)^T\E_{\xi_k}[s_k] &\geq \delta \|\nabla f(x_k)\|_2^2, \\
    \E_{\xi_k}[\|s_k\|_2] &\leq \eta \|\nabla f(x_k)\|_2,\ \ \text{and} \\
    \E_{\xi_k}[\|s_k\|_2^2] &\leq M_{s1} + M_{s2}\|\nabla f(x_k)\|_2^2 \label{eq.stochastic.s.2nd}
  \end{align}
\esubequations
for some $\delta \in (0,1]$, $\eta \in [1,\infty)$, and $(M_{s1},M_{s2}) \in (0,\infty) \times (1,\infty)$ that are all independent of $k$.  For example, as in a stochastic gradient method, these conditions are satisfied if $s_k$ is an unbiased estimate of $\nabla f(x_k)$ with second moment bounded as in~\eqref{eq.stochastic.s.2nd}.  They are also satisfied in the context of a stochastic Newton method wherein a stochastic gradient estimate is multiplied by a stochastic inverse Hessian estimate, assuming that the latter is conditionally uncorrelated with the former and has eigenvalues contained within an interval of the positive real line uniformly over all $k \in \N{}_+$.

Let us also define $\xi_k^H$, conditionally uncorrelated with $\xi_k$ given $x_k$, as a random variable representing a seed for generating an unbiased Hessian estimate $H_k$ such that $\E_{\xi_k^H}[H_k] = \nabla^2 f(x_k)$.  We use~$H_k$ to compute a direction $d_k$.  For the purpose of ideally following a direction of negative curvature (for the true Hessian), we ask that $d_k$ satisfies a curvature condition similar to that used in the deterministic setting.  Importantly, however, the second moment of $d_k$ must be bounded similar to that of $s_k$ above.  Overall, with~$\lambda_k$ being the left-most eigenvalue of $H_k$, we set $d_k \gets 0$ if $\lambda_k \geq 0$, and otherwise require the direction $d_k$ to satisfy
\bsubequations\label{def.stochastic.d}
  \begin{align}
    d_k^TH_kd_k &\leq \gamma\lambda_k\|d_k\|_2^2 < 0\ \ \text{given $H_k$} \label{def.stochastic.d.curv} \\
    \text{and}\ \ \E_{\xi_k^H}[\|d_k\|_2^2] &\leq M_{d1} + M_{d2}\|\nabla f(x_k)\|_2^2
  \end{align}
\esubequations
for some $\gamma \in (0,1]$ and $(M_{d1},M_{d2}) \in (0,\infty) \times (1,\infty)$ that are independent of $k$.  One manner in which these conditions can be satisfied is to compute $d_k$ as an eigenvector corresponding to the left-most eigenvalue $\lambda_k$, scaled such that $\|d_k\|_2^2$ is bounded above in proportion to the squared norm $\|s_k\|_2^2$ where $s_k$ satisfies \eqref{def.stochastic.s}.

Note that our conditions in~\eqref{def.stochastic.d} do \emph{not} involve expected descent with respect to the true gradient at $x_k$.  This can be viewed in contrast to \eqref{def.d}, which involves~\eqref{def.d.descent}.  The reason for this is that, in a practical setting, such a condition might not be verifiable without computing the exact gradient explicitly, which might be intractable or prohibitively expensive.  Instead, without this restriction, a critical component of our algorithm is the generation of an independent random scalar~$\omega_k$ uniformly distributed in $[-1,1]$.  With this choice, one finds $\E_{(\xi_k^H,\omega_k)}[\omega_kd_k] = 0$ such that, effectively, $\omega_kd_k$ merely adds noise to $s_k$ in a manner that ideally follows a negative curvature direction.  This leads to Algorithm~\ref{alg.stochastic.2step} below.

\balgorithm[ht]
  \caption{Two-Step Method for Stochastic Optimization}
  \label{alg.stochastic.2step}
  \balgorithmic[1]
    \Require $x_1 \in \R{n}$, $\{\alpha_k\} \subset (0,\infty)$, and $\{\beta_k\} \subset (0,\infty)$
    \For{\textbf{all} $k \in \N{}_+$}
      \State generate uncorrelated random seeds $\xi_k$ and $\xi_k^H$
      \State generate $\omega_k$ uniformly in $[-1,1]$
      \State set $s_k$ satisfying~\eqref{def.stochastic.s}
      \State set $d_k$ satisfying~\eqref{def.stochastic.d}
      \State set $x_{k+1} \gets x_k + \alpha_k s_k + \beta_k \omega_k d_k$ \label{step.stochastic.update}
     \EndFor
  \ealgorithmic
\ealgorithm

Algorithm~\ref{alg.stochastic.2step} maintains the convergence guarantees of a standard stochastic gradient (SG) method.  To show this, let us first prove the following lemma.

\blemma
  For all $k \in \N{}_+$, it follows that
  \bequation\label{eq.fundamental_inequality}
    \baligned
        &\ \E_{(\xi_k,\xi_k^H,\omega_k)}[f(x_{k+1})] - f(x_k) \\
    \leq&\ -(\delta\alpha_k - \thalf L M_{s2}\alpha_k^2 - \tfrac16 L M_{d2}\beta_k^2) \|\nabla f(x_k)\|_2^2 + \thalf L M_{s1} \alpha_k^2 + \tfrac16 L M_{d1} \beta_k^2.
    \ealigned
  \eequation
\elemma
\bproof
  From Lipschitz continuity of $\nabla f$, it follows that
  \bequationNN
    \baligned
      f(x_{k+1}) - f(x_k)
        &= f(x_k + \alpha_k s_k + \beta_k \omega_k d_k) - f(x_k) \\
        &\leq \nabla f(x_k)^T(\alpha_k s_k + \beta_k \omega_k d_k) + \thalf L\|\alpha_k s_k + \beta_k \omega_k d_k\|_2^2.
    \ealigned
  \eequationNN
  Taking expectations with respect to the distribution of the random quantities $(\xi_k,\xi_k^H,\omega_k)$ given $x_k$ and using \eqref{def.stochastic.s}--\eqref{def.stochastic.d}, it follows that
  \bequationNN
    \baligned
        &\ \E_{(\xi_k,\xi_k^H,\omega_k)}[f(x_{k+1})] - f(x_k) \\
    \leq&\ \nabla f(x_k)^T(\alpha_k\E_{\xi_k}[s_k] + \beta_k\E_{(\xi_k^H,\omega_k)}[\omega_kd_k]) \\
        &\ + \thalf L(\alpha_k^2\E_{\xi_k}[\|s_k\|_2^2] + \beta_k^2\E_{(\xi_k^H,\omega_k)}[\|\omega_kd_k\|_2^2]) \\
        &\ + L \alpha_k\beta_k\E_{(\xi_k,\xi_k^H,\omega_k)}[s_k^T(\omega_k d_k)] \\
    \leq&\ -\alpha_k\delta\|\nabla f(x_k)\|_2^2 \\
        &\ + \thalf L \alpha_k^2 (M_{s1} + M_{s2}\|\nabla f(x_k)\|_2^2) + \tfrac16 L \beta_k^2 (M_{d1} + M_{d2}\|\nabla f(x_k)\|_2^2).
    \ealigned
  \eequationNN
  Rearranging, we reach the desired conclusion. \qed
\eproof

From this lemma, we obtain a critical bound similar to one that can be shown for a generic SG method; e.g., see Lemma~4.4 in~\cite{BottCurtNoce16}.  Hence, following analyses such as that in \cite{BottCurtNoce16}, one can show that the \emph{total} expectation for the gap between~$f(x_k)$ and a lower bound for $f$ decreases.  For instance, we prove the following theorem using similar proof techniques as for Theorems~4.8 and 4.9 in \cite{BottCurtNoce16}.

\btheorem
  Suppose that Algorithm~\ref{alg.stochastic.2step} is run with $\alpha_k = \beta_k = \bar\alpha$ for all $k \in \N{}_+$ where
  \bequation\label{eq.alphabar}
    0 < \bar\alpha \leq \frac{\delta}{2L\max\{M_{s2},M_{d2}\}}.
  \eequation
  Then, for all $K \in \N{}_+$, one has that
  \bequationNN
    \baligned
      \E\left[ \frac{1}{K} \sum_{k=1}^K \|\nabla f(x_k)\|_2^2 \right] \leq&\ \frac{2\bar\alpha L \max\{M_{s1},M_{d1}\}}{\delta} + \frac{2(f(x_1) - f_{\inf})}{K\delta\bar\alpha} \\
      \overset{K\to\infty}{\xrightarrow{\hspace{20pt}}}&\ \frac{2\bar\alpha L \max\{M_{s1},M_{d1}\}}{\delta}.
    \ealigned
  \eequationNN
  On the other hand, if Algorithm~\ref{alg.stochastic.2step} is run with $\{\alpha_k\}$ satisfying
  \bequation\label{eq.alpha_special}
    \sum_{k=1}^\infty \alpha_k = \infty\ \ \text{and}\ \ \sum_{k=1}^\infty \alpha_k^2 < \infty
  \eequation
  and $\{\beta_k\} = \{\chi\alpha_k\}$ for some $\chi \in (0,\infty)$, then it holds that
  \bequationNN
    \baligned
      \lim_{K\to\infty} \E\left[ \sum_{k=1}^K \alpha_k \|\nabla f(x_k)\|_2^2 \right] &< \infty,
    \ealigned
  \eequationNN
  from which it follows, with $A_K := \sum_{k=1}^K \alpha_k$, that
  \bequationNN
     \lim_{K\to\infty} \E\left[ \frac{1}{A_K} \sum_{k=1}^K \alpha_k \|\nabla f(x_k)\|_2^2 \right] = 0, 
  \eequationNN
  which implies that $\liminf_{k\to\infty} \E[\|\nabla f(x_k)\|_2^2] = 0$.
\etheorem
\bproof
  First, suppose that Algorithm~\ref{alg.stochastic.2step} is run with $\alpha_k = \beta_k = \bar\alpha$ for all $k \in \N{}_+$.  Then, taking total expectation in \eqref{eq.fundamental_inequality}, it follows with \eqref{eq.alphabar} and since $\thalf > \tfrac16$ that
  \bequationNN
    \baligned
      &\ \E[f(x_{k+1})] - \E[f(x_k)] \\
      \leq&\ -(\delta - L \max\{M_{s2},M_{d2}\} \bar\alpha) \bar\alpha \E[\|\nabla f(x_k)\|_2^2] + L \max\{M_{s1},M_{d1}\} \bar\alpha^2 \\
      \leq&\ -\thalf\delta\bar\alpha \E[\|\nabla f(x_k)\|_2^2] + L \max\{M_{s1},M_{d1}\} \bar\alpha^2.
    \ealigned
  \eequationNN
  Summing both sides of this inequality for $k \in \{1,\dots,K\}$ yields
  \bequationNN
    \baligned
      f_{\inf} - f(x_1) &\leq \E[f(x_{K+1})] - f(x_1) \\
      &\leq -\thalf\delta\bar\alpha \sum_{k=1}^K \E[\|\nabla f(x_k)\|_2^2] + K L \max\{M_{s1},M_{d1}\} \bar\alpha^2.
    \ealigned
  \eequationNN
  Rearranging and dividing further by $K$ yields the desired conclusion.
  
  Now suppose that Algorithm~\ref{alg.stochastic.2step} is run with $\{\alpha_k\}$ and $\{\beta_k\} = \{\chi\alpha_k\}$ such that the former satisfies~\eqref{eq.alpha_special}.  Since the second condition in \eqref{eq.alpha_special} ensures that $\{\alpha_k\} \searrow 0$, we may assume without loss of generality that, for all $k \in \N{}_+$,
  \bequationNN
    \alpha_k \leq \min\left\{\frac{\delta}{2LM_{s2}},\frac{3\delta}{2LM_{d2}\chi^2}\right\}.
  \eequationNN
  Thus, taking total expectation in \eqref{eq.fundamental_inequality} leads to
  \bequationNN
    \baligned
      &\ \E[f(x_{k+1})] - \E[f(x_k)] \\
      \leq&\ -(\delta\alpha_k - \thalf L M_{s2}\alpha_k^2 - \tfrac16 L M_{d2}\beta_k^2) \E[\|\nabla f(x_k)\|_2^2] + \thalf L M_{s1} \alpha_k^2 + \tfrac16 L M_{d1} \beta_k^2 \\
      \leq&\ -\thalf\delta\alpha_k \E[\|\nabla f(x_k)\|_2^2] + (\thalf L M_{s1} + \tfrac16 L M_{d1} \chi^2) \alpha_k^2.
    \ealigned
  \eequationNN
  Summing both sides for $k \in \{1,\dots,K\}$ yields
  \bequationNN
    \baligned
      f_{\inf} - f(x_1)
        &\leq \E[f(x_{K+1})] - f(x_1) \\
        &\leq -\thalf\delta \sum_{k=1}^K \alpha_k \E[\|\nabla f(x_k)\|_2^2] + (\thalf L M_{s1} + \tfrac16 L M_{d1} \chi^2) \sum_{k=1}^K \alpha_k^2,
    \ealigned
  \eequationNN
  from which it follows that
  \bequationNN
    \sum_{k=1}^K \alpha_k \E[\|\nabla f(x_k)\|_2^2] \leq \frac{2(f(x_1) - f_{\inf})}{\delta} + \frac{LM_{s1} + \tfrac13 L M_{d1} \chi^2}{\delta} \sum_{k=1}^K \alpha_k^2.
  \eequationNN
  The second of the conditions in \eqref{eq.alpha_special} implies that the right-hand side here converges to a finite limit when $K \to \infty$.  Then, the rest of the desired conclusion follows since the first of the conditions in \eqref{eq.alpha_special} ensures that $A_K \to \infty$ as $K \to \infty$. \qed
\eproof

%************
% Subsection
%************
\subsection{Dynamic Method}\label{sec.stochastic.dynamic}

Borrowing ideas from the two-step deterministic method in \S\ref{sec.deterministic.two-step}, the dynamic deterministic method in~\S\ref{sec.deterministic.dynamic}, and the two-step stochastic method in~\S\ref{sec.stochastic.two-step}, we propose the dynamic method for stochastic optimization presented as Algorithm~\ref{alg.stochastic.dynamic} below.  After computing a stochastic Hessian estimate $H_k \in \R{n \times n}$ and an independent stochastic gradient estimate $g_k \in \R{n}$, the algorithm employs the conjugate gradient (CG) method \cite{NoceWrig06} for solving the linear system $H_ks = -g_k$ with the starting point $s \gets 0$.  If $H_k \succ 0$, then it is well known that this method solves this system in at most $n$ iterations (in exact arithmetic).  However, if $H_k \not\succ 0$, then the method might encounter a direction of nonpositive curvature, say $p \in \R{n}$, such that $p^TH_kp \leq 0$.  If this occurs, then we terminate CG immediately and set $d_k \gets p$.  This choice is made rather than spend any extra computational effort attempting to approximate an eigenvector corresponding to the left-most eigenvalue of $H_k$.  Otherwise, if no such direction of nonpositive curvature is encountered, then the algorithm chooses $d_k \gets 0$.  In either case, the algorithm sets $s_k$ as the final CG iterate computed prior to termination.  (The only special case is when nonpositive curvature is encountered in the first iteration; then, $s_k \gets -g_k$ and $d_k \gets 0$.)

For setting values for the dynamic parameters $\{L_k\}$ and $\{\sigma_k\}$, which in turn determine the stepsize sequences $\{\alpha_k\}$ and $\{\beta_k\}$, the algorithm employs stochastic function value estimates.  In this manner, the algorithm can avoid computing the exact objective value at any point (which is often not tractable).  In the statement of the algorithm, we use $f_k : \R{n} \to \R{}$ to indicate a function that yields a stochastic function estimate during iteration $k \in \N{}_+$.

\balgorithm[ht]
  \caption{Dynamic Method for Stochastic Optimization}
  \label{alg.stochastic.dynamic}
  \balgorithmic[1]
    \Require $x_1 \in \R{n}$ and $(L_1,\sigma_1) \in (0,\infty) \times (0,\infty)$
    \For{\textbf{all} $k \in \N{}_+$}
        \State generate a stochastic gradient $g_k$ and stochastic Hessian $H_k$
        
        \State run CG on $H_ks = -g_k$ to compute $s_k$ and $d_k$ (as described in the text above)
        \State set $\alpha_k \gets 1/L_k$ and $\beta_k \gets 1/\sigma_k$
        \State set $\xhat_k \gets x_k + \alpha_k s_k$
        \If{$f_k(\xhat_k) > f_k(x_k)$}
          \State set $L_{k+1} \in [L_k,\infty)$ \label{step.Linc}
          \State (optional) reset $\xhat_k \gets x_k$ \label{step.optional_s}
        \Else
          \State set $L_{k+1} \in (0,L_k]$
        \EndIf
        \State set $x_{k+1} \gets \xhat_k + \beta_k d_k$
        \If{$f_k(x_{k+1}) > f_k(\xhat_k)$}
          \State set $\sigma_{k+1} \in [\sigma_k,\infty)$ \label{step.sigmainc}
          \State (optional) reset $x_{k+1} \gets \xhat_k$ \label{step.optional_d}
        \Else
          \State set $\sigma_{k+1} \in (0,\sigma_k]$
        \EndIf
     \EndFor  
  \ealgorithmic
\ealgorithm

As previously mentioned, we do not claim convergence guarantees for this method.  However, we believe that it is well motivated by our previous algorithms.  One should also note that the per-iteration cost between this method and an inexact Newton-CG method is negligible since any computed $d_k \neq 0$ comes essentially for free from the CG routine.  The only significant extra cost might come from the stochastic function estimates, though these can be made cheaper than any stochastic gradient estimate and might even be ignored completely if one is able to tune fixed values for $L$ and $\sigma$ that work well for a particular application.

%************
% Subsection
%************
\subsection{Numerical Experiments}\label{sec.stochastic.numerical}

We implemented Algorithm~\ref{alg.stochastic.dynamic} in Python 2.7.13.  As a test problem, we trained a convolutional neural network to classify handwritten digits in the well-known \texttt{mnist} dataset; see~\cite{LeCuBottBengHaff98}.  Our neural network, implemented using \texttt{tensorflow},\footnote{\href{https://www.tensorflow.org/}{https://www.tensorflow.org/}} is composed of two convolutional layers followed by a fully connected layer.  In each iteration, we computed stochastic gradient and Hessian estimates using independently drawn mini-batches of size 500.  By contrast, the entire training dataset involves 60,000 feature/label pairs.  The testing set involves 10,000 feature/label pairs.

In each iteration, our implementation runs the CG method for at most 10 iterations.  If a direction of nonpositive curvature is found within this number of iterations, then it is employed as $d_k$; otherwise, $d_k \gets 0$.  In preliminary training runs, we occasionally witnessed instances in which a particularly large stochastic gradient led to a large step, which in turn spoiled all previous progress made by the algorithm.  Hence, to control the iterate displacements, we scaled $s_k$ and/or $d_k$ down, when necessary, to ensure that $\|s_k\| \leq 10$ and $\|\beta_kd_k\| \leq 0.2\|\alpha_ks_k\|$.  We witnessed similar poor behavior when the Lipschitz constant estimates were initialized to be too small; hence, for our experiments, we chose $L_1 \gets 80$ and $\sigma_1 \gets 100$ as the smallest values that yielded stable algorithm behavior for both algorithms.  If stochastic function evaluations suggested an objective increase, then an estimate was increased by a factor of 1.2; see Steps~\ref{step.Linc} and \ref{step.sigmainc} in Algorithm~\ref{alg.stochastic.dynamic}.  The implementation never decreases these estimates.  In addition, while the implementation always takes the step $\alpha_ks_k$ (i.e., it does not follow the optional Step~\ref{step.optional_s}), it does reset $x_{k+1} \gets \xhat_k$ (recall Step~\ref{step.optional_d}) if/when the stochastic function estimates predict an increase in $f$ due to the step $\beta_kd_k$.

For comparison purposes, we ran the algorithm twice using the same starting point and initial seeds for the random number generators: once with $\beta_k$ being reset to zero for all $k \in \N{}_+$ (so that no step along $d_k$ is ever taken) and once with it set as in Algorithm~\ref{alg.stochastic.dynamic}.  We refer to the former algorithm as \texttt{SG} since it is a stochastic-gradient-type method that does not explore negative curvature.  We refer to the latter as \texttt{NC} since it attempts to exploit negative curvature.

The training losses as a function of the iteration counter are shown in Figure~\ref{sg_nc_loss}.  As can be seen in the plot, the performance of the two algorithms is initially very similar.  However, after some initial iterations, following the negative curvature steps consistently offers additional progress, allowing \texttt{NC} to reduce the loss and increase both the training and testing accuracy more rapidly than \texttt{SG}.  Eventually, the plots in each of the figures near each other (after around two epochs, when the algorithms were terminated).  This should be expected as both algorithms eventually near stationary points.  However, prior to this point, \texttt{NC} has successfully avoided the early stagnation experienced by \texttt{SG}.

\bfigure[ht]
  \centering
  \includegraphics[width=0.7\textwidth]{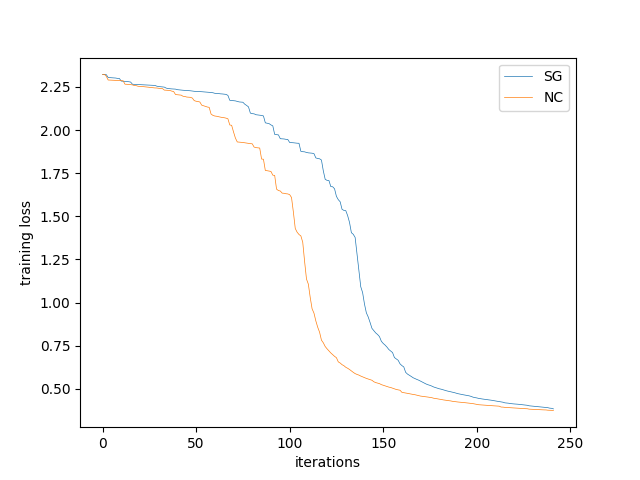}
  \caption{Training loss as a function of the iteration counter when training a convolutional neural network over \texttt{mnist} using an SG-type method (\texttt{SG}) and one that explores negative curvature (\texttt{NC}).}
  \label{sg_nc_loss}
\efigure

\bfigure[ht]
  \centering
  \includegraphics[width=0.48\textwidth]{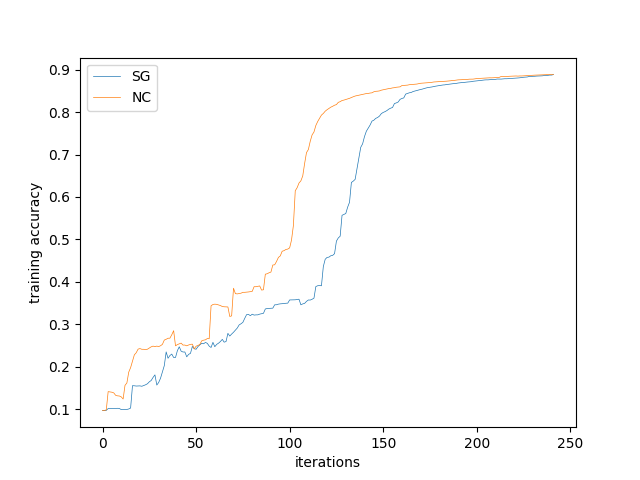}\quad 
  \includegraphics[width=0.48\textwidth]{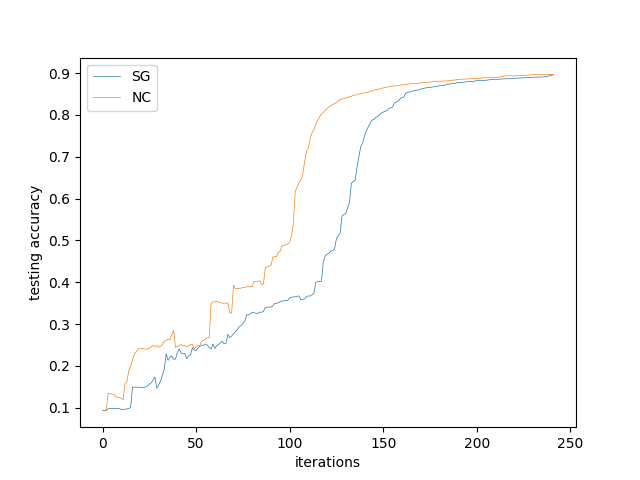}
  \caption{Training (left) and testing (right) accuracy as a function of the iteration counter when training a convolutional neural network over \texttt{mnist} using an SG-type method (\texttt{SG}) and one that explores negative curvature (\texttt{NC}).}
  \label{sg_nc_tr_acc}
\efigure

%*********
% Section
%*********
\section{Conclusion}\label{sec.conclusion}

We have confronted the question of whether it can be beneficial for nonconvex optimization algorithms to compute and explore directions of negative curvature.  We have proposed new algorithmic frameworks based on the idea that an algorithm might alternate or choose between descent and negative curvature steps based on properties of upper-bounding models of the objective function.  In the case of deterministic optimization, we have shown that our frameworks possess convergence and competitive complexity guarantees in the pursuit of first- and second-order stationary points, and have demonstrated that instances of our framework outperform descent-step-only methods in terms of finding points with lower objective values typically within fewer iterations and function evaluations.  In the case of stochastic optimization, we have shown that an algorithm that employs ``curvature noise'' can outperform a stochastic-gradient-based approach.

\bibliographystyle{plain}
\bibliography{sg_ncurv}

\end{document}